\documentclass[a4paper,11pt]{article}

\usepackage{amsmath}
\usepackage{amstext}
\usepackage{amsbsy}
\usepackage{amsopn}
\usepackage{amsthm}
\usepackage{amscd}
\usepackage{amsxtra}
\usepackage{amsfonts}

\def\loadpictex%
  {\ifx\grid\undefined \else \let\normalgrid\grid \fi 
   \ifx\axis\undefined \else \let\normalaxis\axis \fi 
   \ifx\undefined\fiverm
     \font\fiverm=cmr5
   \fi
   \ifx\beginpicture\undefined
     \ifx\newenvironment\undefined
       \input pictex.tex     \relax
     \else
       \input prepictex.tex  \relax
       \input pictex.tex     \relax
       \input postpictex.tex \relax
     \fi
   \fi
   \ifx\normalgrid\undefined \else 
     \let\pictexgrid\grid 
     \let\grid\normalgrid
   \fi 
   \ifx\normalaxis\undefined \else 
     \let\pictexaxis\axis 
     \let\axis\normalaxis
   \fi} 

\ifx\eTeXversion\undefined \else \loadpictex \expandafter \endinput \fi 
 
\ifx \undefined \writestatus \input supp-mis.tex \relax \fi

\unprotect

\let\normalnewdimen = \newdimen
\let\normalnewskip  = \newskip

\def\temporarynewdimen {\alloc@1\dimen\dimendef\insc@unt}
\def\temporarynewskip  {\alloc@2\skip \skipdef \insc@unt}

\temporarynewdimen\!dimenA       
\temporarynewdimen\!dimenB       
\temporarynewdimen\!dimenC       
\temporarynewdimen\!dimenD       
\temporarynewdimen\!dimenE       
\temporarynewdimen\!dimenF       
\temporarynewdimen\!dimenG       
\temporarynewdimen\!dimenH       
\temporarynewdimen\!dimenI       
\temporarynewdimen\!dxpos        
\temporarynewdimen\!dypos        
\temporarynewdimen\!xloc         
\temporarynewdimen\!xpos         
\temporarynewdimen\!yloc         
\temporarynewdimen\!ypos         
\temporarynewdimen\!zpt          

\temporarynewdimen\linethickness 

\def\newdimen#1%
  {\ifx#1\undefined
     \ifnum\count11>\count12\relax
       \temporarynewskip#1\relax
     \else
       \temporarynewdimen#1\relax
     \fi
   \else
   \fi}

\def\dimeninput#1 %
  {\message{[before: d=\the\count11,s=\the\count12]}%
   \input #1 \relax
   \message{[after: d=\the\count11,s=\the\count12]}}%

\loadpictex

\let\newdimen=\normalnewdimen
\let\newskip =\normalnewskip

\protect 

\usepackage{ifthen}
\usepackage{xspace}
\usepackage{fancyheadings}
\usepackage{layout}
\usepackage{hyperref}
\input xy
\xyoption{matrix}
\xyoption{arrow}
\xyoption{curve}
\newdir{ (}{{}*!/-5pt/@^{(}}

\newcommand{\CZ}{\mathbb{C}}

\newcommand{\QZ}{\mathbb{Q}}
\newcommand{\ZZ}{\mathbb{Z}}

\newcommand{\PZ}{\mathbb{P}}

\newcommand{\GZ}{\mathbb{G}}

\newcommand{\sB}{{\mathcal B}}

\newcommand{\sD}{{\mathcal D}}

\newcommand{\sI}{{\mathcal I}}

\newcommand{\sO}{{\mathcal O}}

\newcommand{\sT}{{\mathcal T}}

\newcommand{\from}{\leftarrow}
\newcommand{\suchthat}{\, | \,}

\newboolean{probleme}
\newcommand{\problem}[1]
           {\ifthenelse{\boolean{probleme}}
                       {{\bf(PROBLEM: #1)\bf}}
                       {}
           }
\newboolean{zukuenftiges}
\newcommand{\zukunft}[1]
           {\ifthenelse{\boolean{zukuenftiges}}
                       {{\bf(AUSBAUM\"OGLICHKEIT: #1)\bf}}
                       {}
           }
\newboolean{extras}
\newcommand{\extra}[1]
           {\ifthenelse{\boolean{extras}}
                       {{\bf EXTRA #1 EXTRA\bf}}
                       {}
           }

\DeclareMathOperator{\Img}{Im}

\DeclareMathOperator{\Hom}{Hom}

\DeclareMathOperator{\rank}{rank}

\DeclareMathOperator{\cliff}{cliff}

\DeclareMathOperator{\GL}{GL}
\DeclareMathOperator{\gl}{gl}
\DeclareMathOperator{\Lie}{Lie}

\DeclareMathOperator{\spin}{spin}
\DeclareMathOperator{\Spin}{Spin}
\DeclareMathOperator{\SO}{SO}
\DeclareMathOperator{\Gr}{Gr}                 

\swapnumbers
\theoremstyle{plain}
\begingroup

\newtheorem{thm}{Theorem}
\newtheorem{cor}[thm]{Corollary}
\newtheorem{lem}[thm]{Lemma}

\newtheorem{prop}[thm]{Proposition}
\newtheorem{conj}[thm]{Conjecture}

\numberwithin{thm}{subsection} 
\endgroup

\newtheorem*{thm*}{Theorem}
\newtheorem*{conj*}{Conjecture}
\newtheorem*{verm*}{Vermutung}

\theoremstyle{definition}
\newtheorem{defn}[thm]{Definition}
\newtheorem{rem}[thm]{Remark}
\newtheorem{example}[thm]{Example}
\newtheorem{notation}[thm]{Notation}

\numberwithin{equation}{section}

\newcommand{\nosubsections}{\renewcommand{\thethm}{\thesection.\arabic{thm}}
                            \setcounter{thm}{0}
                           } 
\newcommand{\yessubsections}{\renewcommand{\thethm}{\thesubsection.\arabic{thm}}
                            \setcounter{thm}{0}
                            } 
\newcommand{\cref}[3]{(\ref{#1}, #2 \ref{#3})}

\date{\today}
\setboolean{probleme}{false}
\begin{document}

\title{Geometric Syzygies of Mukai Varieties and
       General Canonical Curves with Genus $\le 8$}

\author{H.-Chr. Graf v. Bothmer
\thanks{Supported by the Schwerpunktprogramm ``Global Methods in Complex
        Geometry'' of the Deutsche Forschungs Gemeinschaft}}

\maketitle

\tableofcontents

\section{Introduction}
\nosubsections

In this paper we study the syzygies of general canonical curves
$C \subset \PZ^{g-1}$ for $g \le 8$.

In \cite{GL1} Green and Lazarsfeld construct low-rank-syzygies of $C$ from 
special linear systems on $C$. More precisely linear systems
of Clifford index $c$ give a $(g-c-3)$rd syzygy. We call these syzygies
geometric syzygies. Green's conjecture
\cite{GreenKoszul} paraphrased in this way is

\begin{conj*}[Green]
Let $C$ be a canonical curve, then 
\[
\text{$C$ has no geometric $k$th syzygies}
\iff
\text{$C$ has no $k$th syzygies at all}
\]
\end{conj*}

This conjecture as received a lot of attention in the last years
and it is now known in many cases \cite{Petri}, \cite{GreenKoszul}, \cite{Sch86},
\cite{Sch88}, \cite{Voisin}, \cite{Sch91}, \cite{Ehbauer}, \cite{HirschR},
\cite{TeixidorGreen}, \cite{VoisinK3}.

A natural generalization of Green's conjecture is 

\begin{conj*}[Geometric Syzygy Conjecture]
Let $C$ be a canonical curve, then the geometric $k$th syzygies span the
space of all $k$th syzygies.
\end{conj*}

Both conjectures are equivalent for $k \ge \frac{g-3}{2}$ since a 
general canonical curve has no linear systems of Clifford index
$c \le \frac{g-3}{2}$.

The geometric syzygy conjecture is therefore true, where Green's conjecture
is known. Furthermore the case $k=0$ (geometric quadrics) was proved
by \cite{AM} for general canonical curves, and by \cite{GreenQuadrics}
for all canonical curves. The case $k=1$ was done for general canonical curves
of genus $g \ge 9$ in \cite{HC}.

In this paper we attack the cases $k=1, g=6,7$ and $k=2, g=8$ for
general canonical curves.

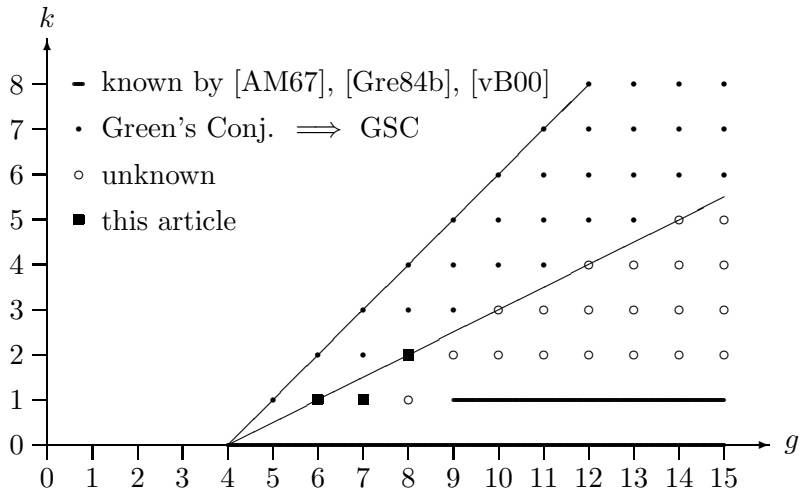
\begin{center}
\mbox{
  \beginpicture   
\setlinear 

\setcoordinatesystem units <6mm,6mm> point at -150 0.5
\unitlength6mm
\setplotarea x from 0 to 16, y from 0 to 9
\axis left 
   ticks numbered from 0 to 8 by 1 /
\put{\vector(0,1){9}} [Bl] at 0 0
\put{$k$} at 0 9.5
\axis bottom
   ticks numbered from 0 to 15 by 1 /
\put {\vector(1,0){16}} [Bl] at 0 0 
\put{$g$} at 16.5 0

\put{\line(1,1){8}} [Bl] at 4 0
\put{\line(2,1){11}} [Bl] at 4 0

\setplotsymbol ({\circle{0.05}} [Bl])
\setsolid
\plot 4 0  8 0  15 0 /
\plot 9 1  10 1  15 1 /

\putrectangle corners at 5.9 0.9 and 6.1 1.1
\putrectangle corners at 6.9 0.9 and 7.1 1.1
\putrectangle corners at 7.9 1.9 and 8.1 2.1

\put{\circle{0.2}} [Bl] at 8 1
\put{\circle{0.2}} [Bl] at 9 2
\put{\circle{0.2}} [Bl] at 10 2
\put{\circle{0.2}} [Bl] at 11 2
\put{\circle{0.2}} [Bl] at 12 2
\put{\circle{0.2}} [Bl] at 13 2
\put{\circle{0.2}} [Bl] at 14 2
\put{\circle{0.2}} [Bl] at 15 2
\put{\circle{0.2}} [Bl] at 10 3
\put{\circle{0.2}} [Bl] at 11 3
\put{\circle{0.2}} [Bl] at 12 3
\put{\circle{0.2}} [Bl] at 13 3
\put{\circle{0.2}} [Bl] at 14 3
\put{\circle{0.2}} [Bl] at 15 3
\put{\circle{0.2}} [Bl] at 12 4
\put{\circle{0.2}} [Bl] at 13 4
\put{\circle{0.2}} [Bl] at 14 4
\put{\circle{0.2}} [Bl] at 15 4
\put{\circle{0.2}} [Bl] at 14 5
\put{\circle{0.2}} [Bl] at 15 5

\put{\circle*{0.1}} [Bl] at 5 1
\put{\circle*{0.1}} [Bl] at 6 2
\put{\circle*{0.1}} [Bl] at 7 2
\put{\circle*{0.1}} [Bl] at 7 3
\put{\circle*{0.1}} [Bl] at 8 3
\put{\circle*{0.1}} [Bl] at 9 3
\put{\circle*{0.1}} [Bl] at 8 4
\put{\circle*{0.1}} [Bl] at 9 4
\put{\circle*{0.1}} [Bl] at 10 4
\put{\circle*{0.1}} [Bl] at 11 4
\put{\circle*{0.1}} [Bl] at 9 5
\put{\circle*{0.1}} [Bl] at 10 5
\put{\circle*{0.1}} [Bl] at 11 5
\put{\circle*{0.1}} [Bl] at 12 5
\put{\circle*{0.1}} [Bl] at 13 5
\put{\circle*{0.1}} [Bl] at 10 6
\put{\circle*{0.1}} [Bl] at 11 6
\put{\circle*{0.1}} [Bl] at 12 6
\put{\circle*{0.1}} [Bl] at 13 6
\put{\circle*{0.1}} [Bl] at 14 6
\put{\circle*{0.1}} [Bl] at 15 6
\put{\circle*{0.1}} [Bl] at 11 7
\put{\circle*{0.1}} [Bl] at 12 7
\put{\circle*{0.1}} [Bl] at 13 7
\put{\circle*{0.1}} [Bl] at 14 7
\put{\circle*{0.1}} [Bl] at 15 7
\put{\circle*{0.1}} [Bl] at 12 8
\put{\circle*{0.1}} [Bl] at 13 8
\put{\circle*{0.1}} [Bl] at 14 8
\put{\circle*{0.1}} [Bl] at 15 8

\put{\circle*{0.1}} [Bl] at 0.7 7
\put{\lines{Green's Conj. $\implies$ GSC}} [l] at 1.2 7

\setplotsymbol ({\circle{0.05}} [Bl])
\setsolid
\plot 0.6 8  0.7 8    0.8 8 /
\put{\lines{known by \cite{AM}, \cite{GreenQuadrics}, \cite{HC}}} [l] at 1.2 8

\put{\circle{0.2}} [Bl] at 0.7 6
\put{\lines{unknown}} [l] at 1.2 6

\putrectangle corners at 0.6 4.9 and 0.8 5.1
\put{\lines{this article}} [l] at 1.2 5

\linethickness1.2mm
\putrule from 5.9 1 to 6.1 1
\putrule from 6.9 1 to 7.1 1
\putrule from 7.9 2 to 8.1 2
\putrule from  0.6 5 to  0.8 5

\endpicture   
     }
\end{center}

The starting point of our proof is

\begin{thm*}[Mukai]
Every general canonical curve of genus $7 \le g \le 9$ is a general
linear section of an embedded rational homogeneous variety $M_g$. 
General canonical curves
of genus $6$ are cut out by a general quadric on a general linear 
section of a homogeneous variety $M_6$.
\end{thm*}

Using this we first consider
the schemes of minimal rank first syzygies ($g=6,7$) respectively 
minimal rank second syzygies ($g=8$) of the Mukai varieties $M_g$
using representation theory. It turns out that all these schemes contain
large rational homogeneous varieties.

Passing from Mukai varieties to canonical curves we describe their
schemes of geometric syzygies as determinantal loci on the above
homogeneous varieties.

Using the resolutions of Eagon-Northcott (for $g=6,8$) and Lascoux 
(for $g=7$) we express the cohomology of the corresponding ideal sheafs 
in terms of the cohomology of homogeneous bundles. 
The later cohomology is then calculated with the theorem of Bott. 
 
This calculation shows $h^0(I(1)) = 0$, proving
the geometric syzygy conjecture in these cases. More precisely our
results are:

\theoremstyle{plain}
\newtheorem*{thmmain6*}{\ref{main6} Theorem}
\newtheorem*{thmmain7*}{\ref{main7} Theorem}
\newtheorem*{thmmain8*}{\ref{main8} Theorem}

\begin{thmmain6*} 
The scheme $Z$ of last scrollar syzygies of a general canonical curve
 $C \subset \PZ^5$ of genus $6$ is a configuration of 
$5$ skew lines in $\PZ^4$ that spans the whole $\PZ^{4}$ of
first syzygies of $C$.
\end{thmmain6*}

\newcommand{\spinorsyz}{{S^+_{syz}}}

\begin{thmmain7*}
The scheme $Z$ of last scrollar syzygies of a general canonical curve
 $C \subset \PZ^6$ of genus $7$ is a linearly normal 
ruled surface of degree $84$
on a spinor variety $\spinorsyz \subset \PZ^{15}$. This ruled
surface spans the whole $\PZ^{15}$ of first syzygies of $C$. 
\end{thmmain7*}

\begin{thmmain8*}
The scheme $Z$ of last scrollar syzygies of a general canonical curve
 $C \subset \PZ^7$ of genus $8$ is a configuration of 
$14$ skew conics that lie on a $2$-uple embedded 
$\PZ^5 \hookrightarrow \PZ^{20}$. $Z$ spans the whole $\PZ^{20}$ of
second syzygies of $C$.
\end{thmmain8*}

The paper is organized as follows:

In section \ref{background} we cover some well known background material
on syzygies. In particular we will introduce the rank of a syzygy,
the scheme of minimal rank syzygies and the vector bundle of linear forms in 
subsection \ref{lowrank}. 

In subsection \ref{scrollar} we will show,
that varieties with very low rank syzygies always lie on certain scrolls.
The properties of these scrollar syzygies are studied and their 
connection with Brill-Noether-Theory is explained. 

Subsection
\ref{linearsections} finally describes what happens to a minimal rank
syzygy of a variety $X$ when $X$ is intersected with a general linear subspace.
It will turn out that the rank of $s$ can drop and that this rank can
be calculated by a morphism of vector bundles $\alpha$ involving
the vector bundle of linear forms.

Section \ref{representationtheory} fixes the notations for the use
of representation theory, and the remaining 
three sections treat curves of genus $6$, $7$ and $8$ in
turn. Starting from the respective Mukai varieties
the proof of the geometric syzygy conjecture outlined above is 
given in full detail.

I would like to thank Kristian Ranestad for the many helpful discussions
during my stay at Oslo University. It was there
where most of the ideas of this work where born.

I dedicate this paper to my grandmother Lilly-Maria who introduced me to
mathematics.


\section{Background on Syzygies} 
\label{background}
\yessubsections

In this section we will review some standard facts from the study of
linear syzygies.

\subsection{Syzygies of Low Rank} 
\label{lowrank}

Let $X \subset \PZ(V)$ be a irreducible non degenerate variety,
$I_X$ generated by quadrics and
\[
   \sI_X \from V_0 \otimes \sO(-2) 
         \xleftarrow{\varphi_1} V_1 \otimes \sO(-3)
         \xleftarrow{\varphi_2} \dots
         \xleftarrow{\varphi_m} V_m \otimes \sO(-m-k)
\]
the linear strand of its resolution

\begin{defn}
An element $s \in V_k$ is called a {\sl$k$-th (linear) syzygy of $X$}.
$\PZ(V_k^*)$ is called the {\sl space of $k$-th syzygies}.
\end{defn}

Every linear syzygy $s$ involves a well defined number of linearly independent
linear forms. This number is called the {\sl rank of $s$}. 
In a more formal way we 
have:

\newcommand{\tildephi}{\tilde{\varphi}}
\newcommand{\spann}[1]{\langle #1 \rangle}
\newcommand{\spans}{\spann{s}}

\begin{defn}\label{tildephi}
Let $s \in V_k$ be a syzygy and
\[
   \tildephi_k \colon V_k \to V_{k-1}\otimes V
\]
the map of vector spaces induced by $\varphi_k$. Then the image
of $s$ under this map
\[
   \tildephi_k(s) \in V_{k-1} \otimes V \cong \Hom(V_{k-1}^*,V).
\]
is a map of vector spaces, and its image
\[
    \spans := \Img \tildephi(s) \subset V
\]
is called the {\sl space of linear forms involved in $s$}. Furthermore
\[
    \rank s := \rank \tildephi(s) = \dim \spans
\]
is called the \sl{rank of $s$}. The {\sl zero set $Z_s$ of a syzygy} 
is the linear
space subspace of $\PZ(V)$ where the linear forms involved in $s$ vanish.
\end{defn}

\begin{example}
Consider the rational normal curve $C \subset \PZ^3$ with minimal resolution
\[
    \sI_C 
    \xleftarrow{(yw-z^2,-xw+yz,xz-y^2)} 
    V_0 \otimes \sO(-2)
    \xleftarrow{\left(\begin{smallmatrix}
                      x & y \\
                      y & z \\
                      z & w \\
                \end{smallmatrix}\right)}
    V_1 \otimes \sO(-2)
\]
where $\dim V_0 = 3$ and $\dim V_1=2$.
If $s,t$ is a basis of $V_1$ and $q_1,q_2,q_3$ a basis of $V_0$,
the map $\tildephi_1$ is given by
\[
    \begin{matrix}
       \tildephi_1 \colon & V_1 &\to & V_0 \otimes V \\
                          & s   &\mapsto &  x\otimes q_1 
                                           + y\otimes q_2 
                                           + z\otimes q_3\\ 
                          & t   &\mapsto &  y\otimes q_1 
                                           + z\otimes q_2 
                                           + w\otimes q_3 
    \end{matrix}
\]
With this, the linear forms involved in $s$ are 
$\spans = \langle x,y,z \rangle$, the rank of $s$ is $3$ and the
Zero locus of $s$ is just the single point $Z_s = (0:0:0:1)$.
\end{example}

\newcommand{\ymin}{{Y_{min}}}

To apply geometric methods to the study of low rank syzygies we projectivize
the space of $k$th syzygies $\PZ(V_k^*)$ and give a determinantal description
of the space $\ymin$ of minimal rank syzygies. 
The linear forms involved in these
syzygies defines a vector bundle on $\ymin$:

\begin{defn}\label{mapofvectorbundles}
On the space of $k$th syzygies $\PZ(V_k^*)$ the map of
vector spaces $\tildephi_k$ induces 
a map of vector bundles
\[
    \psi \colon V_{k-1}^* \otimes \sO_{\PZ(V_k^*)} (-1)
                \to V \otimes  \sO_{\PZ(V_k^*)}
\]
that satisfies
\[
       \psi|_s = \tildephi_k(s) \in \Hom(V_{k-1}^*,V)
\]
The determinantal loci $Y_r(\psi) \subset \PZ(V_k^*)$ of $\psi$ are called
{\sl schemes of rank $r$ syzygies}, since the syzygies in their support
have rank $\le r$.

On the {\sl scheme of minimal rank syzygies 
$\ymin := (Y_{r_{min}}(\psi))_{red}$} 
the restricted
map $\psi|_Y$ has constant rank $r_{min}$. Therefore the image
$L := \Img(\psi|_{\ymin})$ is a vector bundle. We call it the 
{\sl vector bundle of linear forms}, since
\[
                 L|_s = \spans \subset V 
\]
for all minimal rank syzygies $s \in \ymin$.
\end{defn} 

\begin{example}
For the rational normal curve $C \subset \PZ^3$ all first syzygies are of rank
$3$ so $\ymin=Y_3=\PZ(V_1^*)=\PZ^1$. The vector bundle of linear forms on this 
$\PZ^1$ is
the image of
\[
      \psi \colon V_0 \otimes \sO_{\PZ^1}(-1) \to V \otimes \sO_{\PZ^1}
\]
where $\psi$ is given by the flipped syzygy matrix
\[
     \begin{pmatrix}
             s & t & 0 & 0 \\
             0 & s & t & 0 \\
             0 & 0 & s & t 
     \end{pmatrix}.
\]
Since $\psi$ is injective we have $L=3 \sO_{\PZ^1}(-1)$.
\end{example}

\begin{defn}
We say that a variety $X$ as above satisfies the 
{\sl $k$th minimal rank conjecture}, if the scheme 
of $k$th minimal rank syzygies 
$\ymin \subset \PZ(V_k^*)$
is non degenerate. 
\end{defn}

This conjecture is easy to verify for certain rational homogeneous varieties:

\begin{prop} \label{minorbit}
Let $X=G/P \subset \PZ(V)$ be a linearly normal homogeneous rational variety.
If the induced representation
\[
        \rho_k \colon G \to GL(V_k^*)
\]
is irreducible then the variety of minimal rank $k$th syzygies $\ymin$ of
$X$ contains the minimal orbit $G/P_k \subset \PZ(V_k^*)$ of this
representation. In particular $X$ satisfies
the $k$-th minimal rank conjecture.
\end{prop}

\begin{proof}
Since the embedding $X \subset \PZ(V)$ is linearly normal, it induces
a representation 
\[
    \rho \colon G \to \GL(V). 
\]
Since $X$ is $G$ invariant, this $G$ also acts on the minimal free
resolution of $\sI_X$ and therefore induces
a representation
\[
      \rho_k \colon G \to \GL(V_k^*)
\]
Now the rank of a syzygy $s \in \PZ(V_k^*)$ is invariant under coordinate 
transformations of $V$ so that the space of minimal syzygies 
$\ymin \subset \PZ(V_k^*)$ is $G$ invariant. It is also 
compact, so it has to contain a minimal orbit $G/P_k$ of $G$ in $\PZ(V_k^*)$.
Since $\rho_k$ is irreducible,there is only
one minimal orbit and this minimal orbit is non degenerate.
\end{proof}

It will turn out later, that the Mukai varieties for $g=6,7,8$ are of this
type.

\subsection{Scrollar Syzygies}
\label{scrollar}

In this section we want to establish a connection between syzygies
of very low rank of a variety $X$, rational scrolls that contain $X$
and pencils of divisors on $X$. 

First we describe how to construct the equations of a rational scroll $S$
containing $X$ from a $k$th syzygy of rank $k+2$. Conversely we
find that the minimal rank $k$th syzygies of this scroll are
of rank $k+2$. In fact there is a $1:1$ correspondence between
the fibers of $S$ and its minimal rank $k$th syzygies.

Secondly we observe that we can construct a scroll $S$ 
containing $X$ from a pencil of divisors on $X$. Conversely the fibers 
of $S$ will cut out a pencil of divisors on $X$.

It will turn out, that in the case minimal degree complete pencils
of a general canonical curve $C \subset \PZ^{g-1}$ all these
correspondences are the same. We will construct an isomorphism 
from the variety of minimal rank syzygies to the corresponding Brill-Noether
locus in this case.

Except for the construction of the above morphism all results of this
section are well known.

Lets start with

\begin{prop} \label{scrollfromsyzygy}
Let $X \subset \PZ(V)$ be a variety as above and $s \in V_k$ a
$k$th syzygy of rank $k+2$. Then there exists a rational
scroll $S \subset \PZ(V)$ of degree $k+2$ and codimension $k+1$
that contains $X$. Furthermore the vanishing set 
$Z_s \subset \PZ(V)$ is a fiber of $S$.
\end{prop}

\begin{proof}
Let $\{x_1,\dots,x_{k-2}\}$ be a basis of $\spans \subset V$. 
Consider the Koszul complex
\[
     \sO(-k-2) \to \spans^* \otimes \sO(-k-1) 
               \to  \dots
               \to \Lambda^{k} \spans^* \otimes \spans^* \sO(-2)
               \to \Lambda^{k+1} \spans^*\otimes \spans^* \sO(-1) 
\]
where the maps are induced by multiplication with the trace
element
\[
        \sum x_i^* \otimes x_i 
\]
Now $s$ induces a map from $\sO(-k-2)$ to $V_k \otimes \sO(-k-2)$ that
lifts to a map of complexes

\begin{center}
\mbox{
\xymatrix{
        \sO(-k-2) \ar[r] \ar[d]
        & \dots \ar[r] 
        & \Lambda^{k} \spans^*\otimes \spans^* \sO(-2) 
          \ar[r]^\alpha \ar[d] \ar[dr]^{\alpha \circ \beta}
        & \Lambda^{k+1} \spans^* \otimes \spans^* \sO(-1) 
          \ar[d]^\beta \\
        V_k \otimes \sO(-k-2)  \ar[r]
        & \dots \ar[r]
        & V_0 \otimes \spans^* \sO(-2) \ar[r]
        & \sO \ar[r]  
        & \sO_X\\
         }
     }
\end{center}

where $\alpha$ is again given by the multiplication with 
$\sum x_i^* \otimes x_i$ and $\beta$ can be described by the multiplication
with $\sum x_i^* \otimes y_i$ where $\{y_1,\dots,v_{k+2}\} \in V$ are
linear forms. We now claim that the image of $\alpha \circ \beta$
is given by the $2\times 2$-minors of
\[
        M = \begin{pmatrix} 
             x_1 & \dots & x_{k-2} \\
             y_1 & \dots & y_{k-2}
            \end{pmatrix}.
\]
More precisely we claim
\[
   (\alpha \circ \beta) (x_{i_1}^* \wedge \dots \wedge x_{i_k}^*) =
   \pm x_1^* \wedge \dots \wedge x_{k+2}^* \otimes
    \det \begin{pmatrix} x_a & x_b \\ y_a & y_b \end{pmatrix}
\]
where $\{i_1,\dots,i_k\} \cup \{a,b\} = \{1,\dots, k+2\}$. 

Without restriction we can assume $(a,b) = (1,2)$ and check
\begin{align*}
      (\alpha \circ \beta) ( x_3^* \wedge \dots x_{k-2}^* )
      &= \beta( (-1)^k 
        x_1^* \wedge x_3^* \wedge \dots \wedge x_{k-2}^* \otimes x_1 +
        x_2^* \wedge x_3^* \wedge \dots \wedge x_{k-2}^* \otimes x_2) \\
      &= (-1)^{2k}         
        x_1^* \wedge \dots \wedge x_{k-2}^* \otimes x_1y_2
        - x_1^* \wedge \dots \wedge x_{k-2}^* \otimes x_2y_1) \\
      &=  x_1^* \wedge \dots \wedge x_{k+2}^* \otimes
        \det \begin{pmatrix} x_1 & x_2 \\ y_1 & y_2 \end{pmatrix}.
\end{align*}
So the $2\times 2$ minors of $M$ are contained in $I_X$. Therefore
the variety $S$ cut out by these minors contains $X$. 

Furthermore $M$ is $1$-generic. Suppose not, then $M$ has after
row and column-operations the form
\[
       M = \begin{pmatrix} x_1 & x_2 & \dots \\
                            0  & y_2 & \dots
           \end{pmatrix}
\]
and $X \subset Z(x_1y_2) \iff X \subset Z(x_1) \text{ or } X \subset Z(y_2)$
since $X$ is irreducible. This is impossible, since $X$ is non degenerate.

Consequently $M$ is $1$-generic and $S$ is a scroll as claimed above.
\end{proof}

Conversely we have

\begin{lem} \label{syzygyfromscroll}
Let $S \subset \PZ(V)$ be a scroll of degree $k+2$ and
codimension $k+1$, $k \ge 1$. Then the minimal rank
$k$th syzygies of $S$ have rank $k+2$ and the space of minimal rank  
$k$th syzygies
\[
      \ymin \subset \PZ(V_k^*)
\]
is isomorphic to the $k$-uple embedding of $\PZ^1$. Furthermore there
is a $1:1$ correspondence between minimal rank syzygies and fibers 
of $S$.
\end{lem}

\begin{proof}
$S$ is cut out by the $2 \times 2$-minors of a $1$-generic matrix
\[
    A = \begin{pmatrix}
         x_{1} & \dots & x_{k+2} \\
         y_{1} & \dots & y_{k+2} 
        \end{pmatrix}
    \quad x_{i},y_i \in V.
\]
Let 
\begin{align*}
      \Phi_A \colon F \otimes G &\to V \\
                    f_i \otimes g_1 &\to x_{i} \\ 
                    f_i \otimes g_2 &\to y_{i} \\ 
\end{align*}
be the corresponding map of vector spaces, where $\dim F = k+2$ and
$\dim G = 2$. Then $I_S$ is resolved by the Eagon-Northcott complex

\[
     I_S \from V_0 \otimes \sO(-2)
         \xleftarrow{\varphi_1} \dots
         \xleftarrow{\varphi_k} V_k \otimes \sO(-k-2) 
         \from 0
\]

with $V_i = \Lambda^{i+2} F \otimes \Lambda^2 G \otimes S^i G$.

On the space of $k$th syzygies 
$\PZ(V_k^*) \cong \PZ(\Lambda^{k+2} F \otimes \Lambda^2 G \otimes S^k G)^*
\cong \PZ^k$ the group $GL(2)$ acts by coordinate transformation of $G$.
The rank of a syzygy is invariant under such operations and therefore
the space of minimal rank syzygies $\ymin \subset \PZ^k$ is invariant
under $GL(2)$. Furthermore $\ymin$ is compact since we are considering
syzygies of minimal rank. Consequently every component of $\ymin$ must
contain the minimal orbit
\[
      G/P \cong \PZ^1 \xrightarrow{\text{$k$-uple}} \PZ^k.
\]
To show $Y=\PZ^1$ we calculate the tangent space of $\ymin$ in 
\[
      s = f_1 \wedge \dots \wedge f_{k+2} 
          \otimes g_1 \wedge g_2
          \otimes (g_1)^k
        =: f \otimes g \otimes (g_1)^k
\]
We recall the determinantal description of $\ymin$. Consider the
map induced by $\varphi_k$
\begin{align*}
     \tildephi \colon V_k &\to V_{k-1} \otimes V \cong \Hom(V_{k-1}^*,V) \\
                f \otimes g \otimes g_1^i g_2^{k-i} 
                &\mapsto
                \sum_{i,j} f_i^*(f) \otimes g \otimes g_j^*(g_1^i g_2^{k-i})
                                              \otimes \Phi(f_i,g_j)
\end{align*}
Since $\Phi(f_i,g_j) \in V$ this induces a map of vector bundles
\[
    \psi \colon V_{k-1}^* \otimes \sO_{\PZ^k}(-1) \to V \otimes \sO_{\PZ^k}
\]
Consider the basis of $V_{k-1}^*$ which is dual to
\[
  \{f_1^*(f)\otimes g \otimes g_1^{k-1},\dots,
  f_{k+2}^*(f) \otimes g \otimes g_1^{k-1},\dots,
  f_{k+2}^*(f) \otimes g \otimes g_2^{k-1}\}.
\]
In this basis $\psi$ is
then given by the matrix
\[
   \left(
   \begin{array}{ccc|ccc}
       s_0 &        &     & s_1 &        &     \\
           & \ddots &     &     & \ddots &     \\
           &        & s_0 &     &        & s_1 \\
       \hline
       s_1 &        &     & s_2 &        &     \\
           & \ddots &     &     & \ddots &     \\
           &        & s_1 &     &        & s_2 \\
       \hline
           &        &     &     &        &     \\
           & \vdots &     &     & \vdots &     \\
           &        &     &     &        &     \\
       \hline
       s_{k-1} &    &     & s_k &        &     \\
           & \ddots &     &     & \ddots &     \\
           &    & s_{k-1} &     &        & s_k 
   \end{array}
   \right)
   \cdot
   \begin{pmatrix}
    x_{1} \\
    \vdots \\
    x_{k+2} \\
    y_{1} \\
    \vdots \\
    y_{k+2}
    \end{pmatrix}
\]
Where $s_i = (f \otimes g \otimes g_1^{k-i}g_2^i)^* \in V_k^*$ is
a basis of $V_k^*$. In these coordinates we have $s=(1:0:\dots:0)$
and therefore 
\[
     \langle s \rangle = \Img(\Psi|_s) = 
     \langle x_{1}, \dots, x_{k+2} \rangle.
\]
Since the matrix $A$ was $1$-generic, the above linear forms are
linearly independent. This shows that $\rank s = k+2$ for
all $s \in \PZ^1 \subset \PZ^k$. 

Now consider the tangent vectors
\[
     s_\epsilon = (1:\epsilon_1:\dots:\epsilon_{k})
\]
Since we can without restriction suppose, that $y_{1}$ is
linearly independent from $x_{1}, \dots, x_{k+2}$, and since
the syzygies of $Y$ all have rank $k+2$, all 
$(k+3) \times (k+3)$-minors of the first $k+3$ columns of
the above matrix have to vanish for $s_\epsilon$:
\[
   \left(
   \begin{array}{ccc|c}
       1   &        &     & \epsilon_1 \\
           & \ddots &     &      \\
           &        & 1   &      \\
       \hline
       \epsilon_1 & &     & \epsilon_2 \\
           & \ddots &     &      \\
           & & \epsilon_1 &   \\
       \hline
           &        &     &   \\
           & \vdots &     &  \vdots  \\
           &        &     &   \\
       \hline
       \epsilon_{k-1} & & & \epsilon_k \\
           & \ddots &     &      \\
       & & \epsilon_{k-1} &   \\
   \end{array}
   \right)
\]
In particular the determinants
\[
    \det
    \left(
    \begin{array}{ccc|c}
       1   &        &     & \epsilon_1 \\
           & \ddots &     &      \\
           &        & 1   &      \\
       \hline
       \epsilon_i & &     & \epsilon_{i+1} \\
    \end{array}
    \right)
    = \epsilon_{i+1}
\]
must vanish, proving $\epsilon_i=0$ for $i\ge2$. Therefore the
tangent space of $\ymin$ in $s$ can be at most $1$-dimensional. By applying
$GL(2)$ we get the same result for all points of $\PZ^1 \subset \ymin$ 
proving $\PZ^1 = \ymin$.

The $1:1$ correspondence is seen as follows.

The fibers of $S$ are the vanishing sets of the generalized
rows
\[
      (\lambda,\mu) \cdot A
\]
of $A$. Consider the syzygy $s=f \otimes g \otimes (\lambda g_1 + \mu g_2)^k \in \PZ^1$.
Then
\[
     \langle s \rangle = \Img(\tildephi(s)) = 
     \langle \lambda v_{11}    + \mu v_{1,2},\dots,
             \lambda v_{k+2,1} + \mu v_{k+2,1} \rangle
     = (\lambda,\mu) \cdot A.
\] 
\end{proof}

The above propositions 
suggest the following definition:

\begin{defn}
Let $X$ be an irreducible, non degenerate variety 
and $I_X$ generated by quadrics. Then 
the $k$th syzygies of rank $k+2$ are {\sl called scrollar syzygies}. 
The schemes $Y_{k+2}$ of these syzygies are 
called {\sl spaces of scrollar syzygies}. 

If $Y_{k+2}$ is nonempty, then this is also the 
space of minimal rank syzygies $\ymin$
since an irreducible, non degenerate variety can not have $k$th syzygies
of rank $k+1$. 
\end{defn}

\begin{rem}
Scrollar syzygies are the easiest example of the geometric 
syzygies constructed by Green and Lazarsfeld in \cite{GL1}.
\end{rem}

We can now make a precise statement of the geometric syzygy conjecture
for general canonical curves. 

\begin{conj}[Geometric Syzygy Conjecture]
Let $C \subset \PZ^{g-1}$ be a general canonical curve of genus
$g$. Then all minimal rank syzygies are scrollar, and the spaces
of scrollar syzygies are non degenerate.
\end{conj}

\begin{rem}
For special canonical curves it is important to consider the non reduced scheme
structure on the space of scrollar syzygies as can be seen in the
case of a curve of genus $6$ with only one $g^1_5$ \cite[p. 174]{AH}.

Also there are geometric $k$th-syzygies in the sense of Green and Lazarsfeld
\cite{GL1} which are not of rank $k+2$. These must also be considered in the
case of special curves. The easiest example of this phenomenon is exhibited 
by the plane quintic curve of genus $6$ \cite{HC}.
\end{rem}

\newcommand{\scrolld}{S_{|D|}}

We now turn to the connection between scrolls and pencils. Here
we restrict ourselves to
the case of a canonical curve $C \subset \PZ^{g-1}$. Let $|D|$
be a complete pencil of degree $d$ on $C$. Then we can consider the
union
\[
    \scrolld = \bigcup_{D'\in |D|} \spann{D'} \subset \PZ^{g-1}
\]
where $\spann{D'}$ is the linear space spanned by $D'$ in $\PZ^{g-1}$.

\begin{prop}
$\scrolld$ is a rational normal scroll of codimension $g-d$ containing $C$.
\end{prop}

\begin{proof}
Since $C$ is canonically embedded, $\PZ^{g-1}=\PZ(H^0(K))$ with
$K$ a canonical divisor on $C$. The set of hyperplanes in $\PZ^{g-1}$
vanishing on $D'$ is therefore $H^0(K-D')=H^1(D')$ and the codimension
of $\spann{D'}$ correspondingly $h^1(D')=g-d+1$ by Riemann-Roch. This is the
same for all $D'\in |D|$. So $\scrolld$ is a rational scroll of codimension
$g-d$. Its equations are given by the $2\times2$-minors of the
$2 \times (g-d+1)$-matrix obtained from the natural map
\[
      H^0(D) \otimes H^0(K-D) \to H^0(K)
\] 
$\scrolld$ contains $C$ since $D$ moves in a pencil.
\end{proof}

Conversely consider a scroll containing $C$. Its fibers cut out
a pencil of divisors on $C$. These pencils are not always
complete:

\begin{prop}\label{scrollpencil}
Let $C \subset \PZ^{g-1}$ be a non hyperelliptic canonical curve
of genus $g$ contained in a scroll $S$ of codimension $c$. Let
$F$ be a fiber of $S$ and $D = C . F$.
Then $|D|$ is a $g^r_d$ with $r\ge 1$ and $d \le g+r-c-1$.
\end{prop}

\begin{proof}
The fibers of $S$ cut out a pencil of divisors linearly equivalent to $D$.
Therefore $r = \dim |D| \ge 1$.

The codimension of a fiber $F$ is $c+1$, so $h^0(K-D) \ge c+1$.
Riemann-Roch now gives
\[
        d = h^0(D)-h^0(K-D) -1 +g \le r+1 -(c+1)-1+g=g+r-c-1
\]
\end{proof}

\newcommand{\laststeptext}{{\lceil \frac{g-5}{2} \rceil}}
\newcommand{\minclifftext}{{\lceil \frac{g-2}{2} \rceil}}
\newcommand{\mindegtext}{{\lceil \frac{g+2}{2}  \rceil}}
\newcommand{\lastaffinefiberdimtext}{{\lceil \frac{g}{2}  \rceil}} 
\newcommand{\lastprojectivefiberdimtext}{{\lceil \frac{g-2}{2}  \rceil}}

\newcommand{\laststepformula}{\Bigl\lceil \frac{g-5}{2} \Bigr\rceil}
\newcommand{\mincliffformula}{\Bigl\lceil \frac{g-2}{2} \Bigr\rceil}
\newcommand{\mindegformula}{\Bigl\lceil \frac{g+2}{2}  \Bigr\rceil}
\newcommand{\lastaffinfiberdimformula}{\Bigl\lceil \frac{g}{2}  \Bigr\rceil}
\newcommand{\lastprojectivefiberdimformula}
                           {\Bigl\lceil \frac{g-2}{2} \Bigr\rceil}

In particular these linear systems have low Clifford index:

\begin{cor} \label{clifford}
In the situation above we have $\cliff(D) \le g-c-2$. If the
corresponding complete linear system $|D|$ is not a pencil, we even have
 $\cliff(D) \le g-c-3$.
\end{cor}

\begin{proof}
\begin{align*}
   \cliff(D) &= d-2r \\
             & \le g+r-c-1-2r \\
             & = g-c-1-r \\ 
             & \le \left\{ 
                   \begin{matrix}
                     g-c-2 & \text{if $r=1$} \\
                     g-c-3 & \text{if $r>1$}
                   \end{matrix}
               \right.
\end{align*}
\end{proof}

\begin{cor}
A general canonical curve $C \in \PZ^{g-1}$ has 
scrollar syzygies only up to step $\laststeptext$. 
\end{cor}

\begin{proof}
Let $s \in V_k$ be a scrollar syzygy in step $k$.
Then the corresponding scroll $S_s$ has codimension  
$k + 1$ by proposition \ref{scrollfromsyzygy}.
The divisor $D$ cut out by a fiber of $S$ has Clifford index 
\[
   \cliff(D) \le  g-k-3
\]
On the other hand it is well known, that on a general canonical curve
all divisors have Clifford index at least $\minclifftext$. Therefore
\[
     k \le g-3-\mincliffformula = \laststepformula
\]
\end{proof}

\begin{cor}
The scrollar syzygy conjecture implies Green's conjecture for
general canonical curves.
\end{cor}

\begin{proof}
Assume the scrollar syzygy conjecture. Then all minimal rank syzygies
are scrollar. But by the corollary above there are no scrollar syzygies
in step $k>\laststeptext$. Therefore there can be no syzygies at all
in theses steps. This is Greens conjecture for the general canonical curve.
\end{proof}

We want now to consider the last step in the resolution of a
general canonical curve, that still allows syzygies:

\begin{defn}
Let $C \subset \PZ^{g-1}$ be a general canonical curve. 
Then the scrollar $\laststeptext$th syzygies of $C$ are called 
the {\sl last scrollar syzygies} of $C$.
\end{defn}

For the last scrollar syzygies everything is as nice as possible.
First we calculate the degree of the corresponding divisors:

\begin{lem} \label{lastdegree}
Let $C \subset \PZ^{g-1}$ be a general canonical curve, a last scrollar 
syzygy, $S$ the corresponding scroll and $D$ the divisor
cut out by the fiber $F_s$ corresponding to $s$. Then $|D|$ is
a complete pencil of degree $\mindegtext$. 
\end{lem}

\begin{proof}
Suppose $|D|$ was not a complete pencil. Then by corollary \ref{clifford}
we would have
\[
   \cliff{D} \le g-\laststepformula-4 = \Bigl \lceil \frac{g-4}{2} 
                                        \Bigr \rceil.
\]
which is impossible for a general canonical curve.
Consequently we have $r=1$ and $\cliff{D} = \minclifftext$ the minimum possible
value. In particular
\[
     d = \cliff(D)+2r = \mindegformula.
\]\end{proof}

This allows us to construct a morphism from the 
space of last scrollar syzygies to the corresponding Brill-Noether-Locus:

\begin{prop}\label{dimension}
Let $C \subset \PZ^{g-1}$ be a general canonical curve, and $\ymin$
its scheme of last scrollar syzygies.
Then there exists an isomorphism
\[
    \zeta \colon \ymin \to C^1_{\mindegtext}
\]
In particular $\ymin$ is a disjoint union of 
\[
     \frac{2}{g+2}{ g \choose \frac{g}{2} }
\]
rational curves if $g$ is odd, and an irreducible ruled surface
over $W^1_{\mindegtext}$ if $g$ is odd.
\end{prop}

\begin{proof}
Consider the vector bundle of linear forms $L$ on the variety
of $\ymin$ of last scrollar syzygies. Let $Q$ be the cokernel of
the natural inclusion
\[
     0 \to L \to V\otimes\sO_{\ymin} \to Q \to 0.
\]
$Q$ is globally generated and has rank $\lastaffinefiberdimtext$. It therefore
induces a morphism
\[
\begin{matrix}
        \alpha \colon &\ymin &\to &G(V,\lastaffinefiberdimtext) \\
                      &   s  &\to &Z_s
\end{matrix}
\]
where $\GZ := G(V,\lastaffinefiberdimtext)$ is the Grassmannian of 
$\lastaffinefiberdimtext$ dimensional quotient spaces of $V$, or 
equivalently the Grassmannian of $\lastprojectivefiberdimtext$ dimensional
linear subspaces of $\PZ^{g-1}$.

Now consider the incidence  Variety
\[
       I = \{ (\PZ^\lastprojectivefiberdimtext,c) 
              \suchthat 
              c \in \PZ^\lastprojectivefiberdimtext \cap C \subset \PZ^{g-1}
           \} 
           \subset \GZ \times C
\]
and the diagram

\begin{center}
\mbox{
\xymatrix{ 
           {\sD} \ar[r] \ar[d]      & I  \ar[r] \ar[d]   & C \\
           {\ymin} \ar[r]^{\alpha}  & {\GZ}                   \\
         }
     }
\end{center}

obtained by base change. $\sD$ is a family of divisors. The fiber over
a scrollar syzygy $s$ is the divisor cut out by the zero locus $Z_s$ of
$s$. Lemma \ref{lastdegree} shows that these divisors all have
degree $d=\mindegtext$ and $r=1$. By the universal property of
$C^r_d$ we obtain a morphism
\[
     \zeta \colon \ymin \to C^1_{\mindegtext}.
\]
To prove the surjectivity of $\zeta$ take let $D \in C^1_\mindegtext$
be a divisor. The scroll $\scrolld$ spanned by $|D|$ has codimension
$g-\mindegtext$ and the fiber $\spann{D}$ corresponds to
a scrollar $\laststeptext$th (last) syzygy $s$ with $Z_s = \spann{D}$ by
lemma \ref{syzygyfromscroll}.
This implies $D \subset Z_s.C$. Equality follows since they
have the same degree by lemma \ref{lastdegree}. 

We are left to prove that $\zeta$ is injective. Assume $s,t$ are
two last scrollar syzygies, whose zero sets $Z_s$ and $Z_t$ cut
out the same divisor $D=Z_s.C=Z_t.C$. Then the scroll $\scrolld$
obtained from the complete pencil $|D|$ is contained in the
scrolls $S_s$ and $S_t$ corresponding to $s$ and $t$. Now all these
scrolls are of the same dimension, so they have to be equal.
Since there is a $1:1$ correspondence between divisors $D'\in |D|$,
fibers $\spans{D'}$ and scrollar syzygies of $\scrolld=S_s=S_t$
we must have $s=t$.

So $\zeta$ is bijective. Now $\ymin$ is reduced by definition, and
since $C$ is a general canonical curve, 
$C^1_\mindegtext$ is normal. So by Zariskis Main Theorem $\zeta$
is an isomorphism.

The description of $\ymin \cong C^1_\mindegtext$ is obtained from 
Brill-Noether-theory.
For a general canonical curve the dimension of $C^1_\mindegtext$
is given by \cite[p. 214]{ACGH}: 
\begin{align*}
   \dim C^1_\mindegtext &= \rho + 1 \\
                        &= g - 2\left(g-\mindegformula+1\right)+1 \\
                        &= 2\mindegformula-1-g \\
                        &= \left\{
                             \begin{matrix}
                              1 & \text{for $g$ even} \\
                              2 & \text{for $g$ odd}
                             \end{matrix}
                          \right.
\end{align*}
Now the Abel-Jacobi map
\[
\begin{matrix}
       \alpha \colon &C^1_\mindegtext &\to &W^1_\mindegtext \\
                     & D              &\mapsto & |D| \\
\end{matrix}
\]
has $\PZ^1$-fibers, so $C^1_\mindegtext$ is a disjoint union
of finitely many $\PZ^1$'s for $g$ even and a ruled surface over
$W^1_\mindegtext$ for $g$ odd. 

In the even case the number of $\PZ^1$'s can be calculated by a formula
of Castelnuovo \cite[p. 211]{ACGH}:
\[
    \deg W^1_\mindegtext = g! \prod_{i=0}^1 \frac{i!}{(g-\mindegtext+1+i)!}
                         =  \frac{2}{g+2}{ g \choose \frac{g}{2} }
\]
\end{proof}


\subsection{Linear Sections}
\label{linearsections}

In this section we want consider general linear sections $X\cap \PZ(W)$
of $X$. It is well known that such a general linear section
has syzygy spaces of the same dimension as $X$. So one can consider
a syzygy $s$ of $X$ also as a syzygy of $X\cap \PZ(W)$. The rank of
this syzygy can change however, if $\PZ(W)$ does not
intersect the zero locus $Z_s$ of $s$ in the expected codimension.

This will lead to a determinantal
description of those minimal rank syzygies of $X$ that drop rank further
when considered as syzygies of $X\cap \PZ(W)$.

If $X$ is a Mukai-Variety and $C=X \cap \PZ^{g-1}$ a general 
canonical curve it will turn out in the remaining sections of this paper,
that these determinantal subvarieties are of expected dimension and
describe the full space of last scrollar syzygies of $C$.

Let now $X \subset \PZ(V)$ be an irreducible, non degenerate variety, $I_X$
generated by quadrics, 
$\PZ(W) \subset \PZ(V)$ a linear subspace and
$\pi \colon V \to W$ the corresponding projection. Consider the intersection
$X \cap \PZ(W)$ and the linear strand
\[
   \sI_{X\cap \PZ(W)/\PZ(W)} 
                          \from W_0 \otimes \sO(-2) 
         \xleftarrow{\varphi_1} W_1 \otimes \sO(-3)
         \xleftarrow{\varphi_2} \dots
         \xleftarrow{\varphi_m} W_m \otimes \sO(-m-k)
\]
of its resolution.

We start by recalling:

\begin{prop}\label{isomorphisms}
The inclusion $\PZ(W) \subset \PZ(V)$ induces linear maps
\[
    \pi_k \colon V_k \to W_k
\]
If $X$ is arithmetically Cohen Macaulay and $X \cap \PZ(W)$ is of
the expected dimension, then all $\pi_k$'s are isomorphisms.
\end{prop} 

\begin{proof} Since $X$ is arithmetically Cohen Macaulay, we have in particular
$H^1(X,\sO_X(q))=0$ for all $q \ge 0$. We can then apply \cite[Thm 3.b.7]{GreenKoszul}
to get the result.
\end{proof}

If we regard a syzygy $s \in \ymin(X)$ as a syzygy of $X \cap \PZ(W)$ in
the way made precise by the preceding proposition, we can calculate
the rank of $s$ there by using the vector bundle of linear forms:

\begin{cor}\label{determinantal}
Let $X$ be ACM, $X \cap \PZ(W)$ of expected dimension, 
$\ymin(X) \subset \PZ(V_k^*) \cong
\PZ(W_k^*)$ the scheme of $k$-th minimal rank syzygies of 
$X$ and $L$ the vector bundle of linear forms on $\ymin(X)$. 
Then there
exists a map of vector bundles
\[
     \alpha \colon L \to W \otimes \sO_{\ymin(X)}
\]
such that the rank of a syzygy $s \in \ymin(X)$ considered
as a syzygy of $X\cap\PZ(W)$ is
\[
     \rank_{X \cap \PZ(W)} s = \rank \alpha|_s
\]
\end{cor}

\begin{proof}
We have a diagram

\begin{center}
\mbox{
\xymatrix{
           V_k \ar[r]^-{\tildephi^V_k} \ar[d]_{\pi_k}    
           & V_{k-1}\otimes V \ar[d]^{\pi_{k-1}\otimes\pi}\\
           W_k \ar[r]^-{\tildephi^W_k}   
           & W_{k-1}\otimes W
         }
     }
\end{center}

Since $\pi_k$ and $\pi_{k-1}$ are isomorphisms, this induces a diagram

\begin{center}
\mbox{
\xymatrix{
           V_{k-1}^* \otimes \sO_{\PZ(V^*_k)}(-1) 
           \ar@{=}[d] \ar[r]^-{\psi^V}
           & V \otimes \sO_{\PZ(V^*_k)} \ar[d]_{\pi} \\
           W_{k-1}^* \otimes \sO_{\PZ(W^*_k)}(-1) \ar[r]^-{\psi^W}
           & W \otimes \sO_{\PZ(W^*_k)} 
         }
     }
\end{center}

on $\ymin$ the top map factors over $L$ yielding

\begin{center}
\mbox{
\xymatrix{
                      &   L \cong \Img(\psi^V) \ar@{ (->}[d] \ar@/^3pc/[dd]^\alpha\\
           V_{k-1}^* \otimes \sO_\ymin(-1) \ar@{->>}[ur]^\beta
           \ar@{=}[d] \ar[r]^-{\psi^V}
           & V \otimes \sO_\ymin \ar[d]_{\pi} \\
           W_{k-1}^* \otimes \sO_\ymin(-1) \ar[r]^-{\psi^W}
           & W \otimes \sO_\ymin 
         }
     }
\end{center}

Since $\beta$ is surjective, we have for every syzygy $s\in \ymin$
\[
     \rank(\psi^W|_s) = \rank(\alpha|_s)
\]
as claimed.
\end{proof}

A stronger statement is true for a general intersection with
quadric hypersurface. Here the dimension of syzygy spaces and
the rank of the syzygies stay the same:

\begin{prop}\label{generalquadric}
Let $X \in \PZ(V)$ be an irreducible variety, $I_X$ generated by quadrics and
$X \cap Q$ the intersection with a general quadric $Q \subset \PZ(V)$. If
\[
    I_X \from V_0 \otimes \sO(-2) 
        \xleftarrow{\varphi_0} V_1 \otimes \sO(-3)
        \from \dots
        \from V_k \otimes \sO(-k-2)
        \from 0
\]
is the linear strand of the resolution of $I_X$, then
\[
    I_{X \cap Q}
        \from (V_0 \oplus Q) \otimes \sO(-2) 
        \xleftarrow{\varphi_0\oplus 0} V_1 \otimes \sO(-3)
        \from \dots
        \from V_k \otimes \sO(-k-2)
        \from 0
\]
is the linear strand of the resolution of $I_{X \cap Q}$. All differentials are
the same except for $\varphi_0$ whose matrix has one more column of zeros.

In particular the $k$th syzygies ($k \ge 1$) of $X \cap Q$ are the same as 
those of $X$ and they have the same ranks.
\end{prop}

\begin{proof}
We prove the proposition on the ring level. Let $R = \CZ[V]$ be the
coordinate ring of $\PZ(V)$ and
\[
   C_\bullet \colon R \from C_1 \from \dots \from 0
\]
the graded minimal free resolution of $R/I_X$. $R/QR$ is resolved by
the Koszul complex
\[
   K_\bullet(Q) \colon R \xleftarrow{Q} R(-2) \from 0
\]
Since $Q$ is a nonzerodivisor in $R/I_X$ the total complex
$C_\bullet(Q) := C_\bullet \otimes K_\bullet(Q)$ is a resolution
of $R/QI_X = R/I_{X\cap Q}$ \cite[Thm 16.4.]{Matsumura}.

For degree reasons the linear strand of $C_\bullet(Q)$ is the same as the one of
$C_\bullet$ except for the first step.
\end{proof}

Later we will apply the last two propositions to Mukai varieties, as
defined by

\begin{thm}[Mukai]
Every general canonical curve of genus $7 \le g \le 9$ is a general
linear section of an embedded rational homogeneous (Mukai) variety $M_g$. 
General canonical curves
of genus $6$ are cut out by a general quadric on a general linear 
section of a homogeneous (Mukai) variety $M_6$.

More explicitly we have

\begin{center}
\begin{tabular}{|c|c|}
\hline
 $g$ & $M_g$  \\
\hline
  $6$ & the Grassmannian $\Gr(2,5) \subset \PZ^9$ \\
  $7$ & the Spinor-Variety $S_{10} \subset \PZ^{15}$ \\
  $8$ & the Grassmannian $\Gr(2,6) \subset \PZ^{14}$ \\
  $9$ & the symplectic Grassmannian $\Gr(3,6,\eta) \subset \PZ^{13}$\\
\hline
\end{tabular}
\end{center}
\end{thm}

\begin{proof} \cite{FanoMukai} \cite{CurvesMukai}
\end{proof}

\section{Representation Theory}
\label{representationtheory}
\nosubsections

As the Mukai varieties are rational homogeneous, the main tool of our
study will we representation theory.
We will use the following notations from Fulton/Harris \cite{Fu} and 
Ottaviani \cite{Ott}:

\newcommand{\lie}[1]{{\mathfrak #1}}
\newcommand{\weyl}{{\mathcal W}}

\begin{notation}
We will denote by

\begin{tabular}{ll}
$G$ & 
a semisimple and simply connected Lie group \\
$P \subset G$ & 
a parabolic subgroup \\
$\lie{p} \subset \lie{g}$ & 
the corresponding Lie algebras\\
$\lie{h} \subset \lie{p} \subset \lie{g}$ & 
a Cartan subalgebra \\
$\{H_i\} \subset \lie{h}$ & 
a basis of $\lie{h}$\\
$\{L_i\} \subset \lie{h}^*$ &
the dual basis to $\{H_i\}$ \\
$R=\{\alpha_i\} \subset \lie{h}^*$ & 
the roots of $\lie{g}$ \\
$R= R^+ \cup R^-$ &
a decomposition into positive and negative roots \\
$\lie{g} = \lie{h} \bigoplus 
\left(\oplus_{\alpha_i\subset R} \lie{g}_{\alpha_i} \right)$ &
the Cartan decomposition \\
$\Delta = \{ \alpha_1, \dots, \alpha_k\} \subset R^+$ &
the set of simple roots \\
$\omega_1,\dots,\omega_k$ &
the corresponding fundamental weights\\ 
$\Sigma \subset \Delta$ &
a subset of simple roots \\
$R^+(\Sigma) = \{ \alpha \in R^+ | 
                  \alpha = \sum_{\alpha_i \not\in \Sigma} p_i \alpha_i \}$ &
the positive roots generated by $\Delta - \Sigma$. \\
$R^-(\Sigma) = \{ \alpha \in R^- | 
                  \alpha = \sum_{\alpha_i \not\in \Sigma} p_i \alpha_i \}$ &
the negative roots generated by $\Delta - \Sigma$. \\
$\lie{p}(\Sigma)$ & 
the subalgebra $\lie{h} \bigoplus 
\left( \oplus_{\alpha \in R^+} \lie{g}_\alpha \right)
\bigoplus 
\left( \oplus_{\alpha \in R^-(\Sigma)} \lie{g}_\alpha \right)$ \\
$P(\Sigma)$ &
the corresponding parabolic subgroup of $G$ \\
$S_{P(\Sigma)}$ &
the semisimple part of $\lie{p}(\Sigma)$ \\
$\weyl$ &
the Weyl group of $\lie{g}$
\end{tabular}
\end{notation}


We also need the notion of a highest weight vector:

\begin{defn}
Let $\rho \colon \lie{g} \to \lie{gl}(V)$ be a representation. A
vector $v \in V$ with
\[
 \rho(\lie{g}_\alpha)(v) = 0 \quad \forall \alpha \in R^+ 
\quad \quad \text{and} \quad \quad
 \rho(\lie{h})(v) = \lambda(H)v, \quad \lambda \in \lie{h}^*
\]
is called a highest weight vector. $\lambda \in \lie{h}^*$ is then
called a highest weight.
\end{defn}

With this we use further notations

\newcommand{\schur}{S}

\newcommand{\youngone}[1]{\ifthenelse{\equal{1}{#1}}{\\\cline{1-1}}{}
                          \ifthenelse{\equal{2}{#1}}{&\\\cline{1-2}}{}
                          \ifthenelse{\equal{3}{#1}}{&&\\\cline{1-3}}{}
                          \ifthenelse{\equal{4}{#1}}{&&&\\\cline{1-4}}{}
                          \ifthenelse{\equal{5}{#1}}{&&&&\\\cline{1-5}}{}
                          \ifthenelse{\equal{6}{#1}}{&&&&&\\\cline{1-6}}{}
                          \ifthenelse{\equal{7}{#1}}{&&&&&&\\\cline{1-7}}{}
                          \ifthenelse{\equal{8}{#1}}{&&&&&&&\\\cline{1-8}}{}
                         }
\newcommand{\younglast}[1]{\ifthenelse{\equal{1}{#1}}{\\\cline{1-1}\end{array}}{}
                           \ifthenelse{\equal{2}{#1}}{&\\\cline{1-2}\end{array}}{}
                           \ifthenelse{\equal{3}{#1}}{&&\\\cline{1-3}\end{array}}{}
                           \ifthenelse{\equal{4}{#1}}{&&&\\\cline{1-4}\end{array}}{}
                           \ifthenelse{\equal{5}{#1}}{&&&&\\\cline{1-5}\end{array}}{}
                           \ifthenelse{\equal{6}{#1}}{&&&&&\\\cline{1-6}\end{array}}{}
                           \ifthenelse{\equal{7}{#1}}{&&&&&&\\\cline{1-7}\end{array}}{}
                           \ifthenelse{\equal{8}{#1}}{&&&&&&&\\\cline{1-8}\end{array}}{}
                          }
\newcommand{\youngstart}[1]{\begin{array}{|c|c|c|c|c|c|c|c|}\cline{1-#1}\youngone{#1}}
\newcommand{\youngonly}[1]{\begin{array}{|c|c|c|c|c|c|c|c|}\cline{1-#1}\younglast{#1}}
\newcommand{\y}[1]{\youngonly{#1}}
\newcommand{\yy}[2]{\youngstart{#1}\younglast{#2}}
\newcommand{\yyy}[3]{\youngstart{#1}\youngone{#2}\younglast{#3}}
\newcommand{\yyyy}[4]{\youngstart{#1}\youngone{#2}\youngone{#3}\younglast{#4}}
\newcommand{\yyyyy}[5]{\youngstart{#1}\youngone{#2}\youngone{#3}\youngone{#4}\younglast{#5}}
\newcommand{\yyyyyy}[6]{\youngstart{#1}\youngone{#2}\youngone{#3}\youngone{#4}\youngone{#5}\younglast{#6}}

\begin{notation}
We denote by

\begin{tabular}{ll}
$\rho_\lambda$ &
the irreducible representation of $\lie{g}$ with highest weight $\lambda$ \\
$\schur_\lambda$ &
the Schur functor for a partition $\lambda$ \\
$\Lambda_\lambda$ &
the Schur functor for the dual partition $\tilde{\lambda}$
\end{tabular}
\end{notation}

\begin{rem}
With this notation we have for example
\[
      \schur_{1,1,1} = \Lambda_{3} = \Lambda^3
\]
where $\Lambda^3$ is the usual exterior product.

Sometimes we represent the Schur functors by young tableaux, 
in this case
\[
        \schur_{1,1,1} = \Lambda_{3} = \yyy{1}{1}{1}
\]
If $\lie{g} = \lie{gl}(V)$ we also have the identity
$
    \rho_{\lambda_1 L_1 + \dots + \lambda_2 L_n}
    = \schur_{\lambda_1,\dots,\lambda_n}
$.
\end{rem}

If $\rho$ is a representation of parabolic subgroup $P \subset G$
with highest weight $\lambda=\lambda_1L_1+\dots+\lambda_nL_n$ 
we use the notation $E_{\rho}$, $E(\lambda)$ or $E(\lambda_1,\dots,\lambda_n)$
for the vector bundle induced by $\rho$ on $G/P$.

\begin{thm}[Matsushima]
A vector bundle E of rank r over $G/P$ is homogeneous if and only if 
there exists a representation $\rho \colon P \to GL(r)$ such
that $E\cong E_\rho$.
\end{thm}

\begin{proof} \cite[Theorem 9.7]{Ott}
\end{proof}

\begin{thm}[Classification of irreducible bundles over $G/P$]
Let $P(\Sigma) \in G$ be a parabolic subgroup and $\omega_1,\dots,\omega_k$ the 
fundamental weights corresponding to the subset of simple roots $\Sigma \subset \Delta$.
Then all irreducible representations of $P(\Sigma)$ are
\[
     V \otimes L^{n_1}_{\omega_1} \otimes \dots \otimes L^{n_k}_{\omega_k} 
\]
where $V$ is a representation of $S_P$ and $n_i \in \ZZ$. $L_{\omega_i}$
are the one dimensional representations of $S_P$ induced by the 
fundamental weights.

The weight lattice of $S_P$ is embedded in the weight lattice of
$G$. If $\lambda$ is the highest weight of $V$, we will
call $\lambda + \sum n_iw_i$ the highest weight of the irreducible
representation of $P(\Sigma)$ above. 
\end{thm}

\begin{proof} \cite[Proposition 10.9 and remark 10.10]{Ott}
\end{proof}

For the cohomology of homogeneous vector bundles we use

\begin{thm}[Bott]
Consider the homogeneous vector bundle $E(\lambda)$ on $X=G/P$ and
$\delta$ the sum of fundamental weights of $G$. Then
\begin{itemize}
\item $H^i(X,E(\lambda))$ vanishes for all $i$ if there
is a root $\alpha$ with $(\alpha, \delta+\lambda)=0$
\item Let $i_0$ be the number of positive roots $\alpha$ with 
$(\alpha, \delta+\lambda) < 0$. Then $H^{i}(X,E(\lambda))$ vanishes
for $i\not=i_0$ and $H^{i_0}(X,E(\lambda))=\rho_{w(\delta+\lambda)-\delta}$
\end{itemize}
where $(.,.)$ denotes the Killing form on $\lie{h}^*$, 
$w(\delta+\lambda)$ is the unique element of the fundamental Weyl chamber 
which is congruent to $\delta+\lambda$ under the action of the Weyl group,
and $\rho_{w(\delta+\lambda)-\delta}$ is the corresponding representation
of $G$.
\end{thm}

\begin{proof}
\cite[Theorem 11.4]{Ott}
\end{proof}

In several proofs of this paper we need to calculate the decomposition of
tensor products of fundamental representations. Formulas for this can be found in 
\cite{fundamentalKempf}. Also we sometimes calculate the dimension of certain
representations. This can be done in various combinatoric ways and by the use
of the Weyl character formula. In this paper all calculations of this kind have 
been checked by the computer program SYMMETRICA \cite{symmetrica} which can be used 
online on the web.

\section{Genus 6}
\yessubsections

\newcommand{\PZdual}{\hat{\PZ}}

The following is well known, but sets the stage for the more involved
computations for Mukai varieties of higher genus.

\subsection{Syzygies of $M_6$}

Let $V$ be a $5$-dimensional vector space with basis $\{v_1,\dots,v_5\}$.
In this section we will abbreviate $\Lambda_\lambda V$ by $\Lambda_\lambda$.

The Mukai variety for genus $6$ is
\[
       M_6 = \Gr(V,2) = \Gr(5,2) \cong \GL(5)/P 
             \subset \PZ(\Lambda_2) \cong \PZ^9
\]
The diagonal matrices $\lie{h}$ form a Cartan subalgebra of $\lie{gl_5}$
and $\lie{p}$. The matrices $H_i= E_{i,i}$ are a basis of $\lie{h}$.
Let $\{L_i\}$ be the dual basis of $\lie{h}^*$.  The positive roots
of $\lie{gl_5}$ are $L_i-L_j$ with $i>j$ and $\omega_i = \sum_{j=1}^i L_j$
are the fundamental weights.

Lets first consider the dimensions of the spaces of linear syzygies. To write
these in compact form we use the MACAULAY-notation \cite{M2}:

\begin{prop} The syzygy-numbers of $M_6 = \Gr(5,2)$
are
\[
\begin{matrix}
 1 & - & - & - \\
 - & 5 & 5 & - \\
 - & - & - & 1 
\end{matrix}
\]
\end{prop}

\begin{proof}
The Grassmannian $\Gr(5,2)$ is cut out by the $4\times4$ Pfaffians of a
generic $5 \times 5$ matrix $A$. It is therefore a codimension $3$ Gorenstein
variety and has the resolution
\[
      I_{\Gr(5,2)} \from V^* \otimes \sO(-2) 
                   \xleftarrow{A} V \otimes \sO(-3)
                   \from \sO(-5)
                   \from 0.
\]
Since $\dim V = 5$ this gives the above syzygy numbers.
\end{proof}

Since $\GL(5)$ acts on $\Gr(5,2)$ we also have an action on the syzygies. 
The corresponding representations are calculated by

\begin{prop} The linear strand of the resolution of
$\Gr(5,2)$ is
\[
   I_{\Gr(5,2)} \from \Lambda_4  \otimes \sO(-2)
              \from \Lambda_{51}  \otimes \sO(-3)
\]
\end{prop}

\begin{proof}
From above we have a linear strand
\[
   I_{\Gr(5,2)} \from V_0 \otimes \sO(-2)
              \from V_1 \otimes \sO(-3)
\]
with $\dim V_0 = \dim V_1 = 5$. $V_0$ is a invariant
(not necessarily irreducible)
subspace of quadrics. This gives
\[
    V_0 \subset S_2(\Lambda_2 ) 
        \subset \Lambda_2  \otimes \Lambda_2 
        =  \Lambda_4  \oplus \Lambda_{31}  \oplus \Lambda_{22}  
\]
where the irreducible components have dimensions $5$, $45$ and $50$
respectively. So $V_0 = \Lambda_4 $. Similarly we have
\[
    V_1 \subset \Lambda_4  \otimes \Lambda_2  
         = \Lambda_{51} 
           \oplus \Lambda_{42} 
\]
where the irreducible components have dimensions $5$ and $45$
respectively. This implies $V_1 = \Lambda_{51}$.
\end{proof}

This allows us to describe the minimal rank syzygies $\ymin$ of $\Gr(5,2)$ 
and the vector bundle of linear forms on $\ymin$

\begin{prop} \label{yminM6}
The scheme of minimal rank first syzygies of $\Gr(5,2)$ is
\[
     \ymin \cong GL(5)/P \cong \PZ(\Lambda_{51}^*) \cong \PZ^4
\]
The bundle of linear forms on $\ymin$ is
\[
     L|_\ymin = E(1,1,0,0,0)^* = \sT_{\PZ^4}(-2) 
\]
$L$ has rank $4$.
\end{prop}

\begin{proof}
From proposition \ref{minorbit} we know, that 
$\ymin \subset \PZ(\Lambda_{51}^*) \cong \PZ^4$ must contain the minimal
orbit of $GL(5)$ in $\PZ(\Lambda_{51}^*)$ under the action
\[
       \rho \colon GL(5) \to GL(\Lambda_{51}^*).
\] 
Here this orbit $GL(5)/P$ is the
whole $\PZ^4$ such that $\ymin=\PZ^4$.

To describe the vector bundle of linear forms on $\ymin=GL(5)/P$ we have to
determine the action of $P$ on a fiber of $L$. We start by considering
the dual actions $\rho^*$ of $GL(5)$ and $P$
on $\Lambda_{51}$. The parabolic subgroup $P$ in its standard representation
is then the set of
matrices that fix a given $\PZ^0$, i.e matrices of the
form
\[
   \begin{pmatrix}
         * & * & * & * & * \\
         0 & * & * & * & * \\
         0 & * & * & * & * \\
         0 & * & * & * & * \\
         0 & * & * & * & * 
   \end{pmatrix}
\]
The semisimple part of $P$ is $S_P = \GL(1) \times \GL(4)$
where $\GL(1)$ acts on $\langle v_1 \rangle$ and $\GL(4)$ acts on
$\langle v_2,v_3,v_4,v_5 \rangle$ in the standard way.

The maximal weight vector
\[
   s = \begin{tabular}{|c|c|}\hline
   1 & 1 \\ \hline
   2 \\ \cline{1-1}
   3 \\ \cline{1-1}
   4 \\ \cline{1-1}
   5 \\ \cline{1-1}
   \end{tabular} 
   \in \Lambda_{51}
\]
is a syzygy in $\ymin$. To determine the fiber $L|_s$ of $L$ over
$s$ we restrict the map
\[
    \psi \colon \Lambda_4^* \otimes \sO_{\PZ(\Lambda_{51}^*)}(-1)
                \to \Lambda_2 \otimes \sO_{\PZ(\Lambda_{51}^*)}
\]
from definition \ref{mapofvectorbundles} to $s$. This gives
\[
             \psi|_s =\tildephi(s) \in \Hom(\Lambda_4^*,\Lambda_2) 
                                       \cong \Lambda_4 \otimes \Lambda_2
\]
where
\[
             \tildephi \colon 
             \Lambda_{51} 
             \hookrightarrow 
             \Lambda_4 \otimes \Lambda_2.
\]
Using Young-diagrams we get

\begin{align*}
   \tildephi \colon \yyyyy{2}{1}{1}{1}{1} 
   &\hookrightarrow  
   \yyyy{1}{1}{1}{1}\otimes \yy{1}{1}
   \\
   \begin{tabular}{|c|c|}\hline
   1 & 1 \\ \hline
   2 \\ \cline{1-1}
   3 \\ \cline{1-1}
   4 \\ \cline{1-1}
   5 \\ \cline{1-1}
   \end{tabular} 
   &\mapsto
   \begin{tabular}{|c|}\hline
   1 \\ \hline
   2 \\ \hline
   3 \\ \hline
   4 \\ \hline
   \end{tabular} 
   \otimes
   \begin{tabular}{|c|}\hline
   1 \\ \hline
   5 \\ \hline
   \end{tabular} 
   -
   \begin{tabular}{|c|}\hline
   1 \\ \hline
   2 \\ \hline
   3 \\ \hline
   5 \\ \hline
   \end{tabular} 
   \otimes
   \begin{tabular}{|c|}\hline
   1 \\ \hline
   4 \\ \hline
   \end{tabular} 
   +
   \begin{tabular}{|c|}\hline
   1 \\ \hline
   2 \\ \hline
   4 \\ \hline
   5 \\ \hline
   \end{tabular} 
   \otimes
   \begin{tabular}{|c|}\hline
   1 \\ \hline
   3 \\ \hline
   \end{tabular} 
   -
   \begin{tabular}{|c|}\hline
   1 \\ \hline
   3 \\ \hline
   4 \\ \hline
   5 \\ \hline
   \end{tabular} 
   \otimes
   \begin{tabular}{|c|}\hline
   1 \\ \hline
   2 \\ \hline
   \end{tabular} 
\end{align*}

Consequently the fiber $L|_s$ of the line bundle of linear forms is
\[
    L|_s = \Img \tildephi(s) = 
   \langle 
   \begin{tabular}{|c|}\hline
   1 \\ \hline
   5 \\ \hline
   \end{tabular}, 
   \begin{tabular}{|c|}\hline
   1 \\ \hline
   4 \\ \hline
   \end{tabular},
   \begin{tabular}{|c|}\hline
   1 \\ \hline
   3 \\ \hline
   \end{tabular}, 
   \begin{tabular}{|c|}\hline
   1 \\ \hline
   2 \\ \hline
   \end{tabular}
   \rangle
   = \langle
   v_1 \wedge v_5, v_1 \wedge v_4, v_1 \wedge v_3, v_1 \wedge v_2
   \rangle
\]
and $L$ is of rank $4$. $S_P$ acts irreducibly on this fiber, and
$v_1 \wedge v_2$ is the maximal weight vector of weight $L_1 + L_2$.
(Notice that the weights with respect to $G$ are the same
as the ones with respect to $S_P$ since we can use the same
Cartan subalgebra $\lie{h} \subset Lie S_P \subset \lie{p} \subset \lie{g}$).

Therefore
\[
          E_{\rho^*} = E(1,1,0,0,0)
\]
and
\[
          L = E_\rho = E(1,1,0,0,0)^*
\]
Now the tautological sequence of $\PZ^4$ is 
\[
        0 \to E(0,1,0,0,0) 
          \to V \otimes \sO_{\PZ^4} 
          \to E(1,0,0,0,0) 
          \to 0
\]
i.e. $E(1,0,0,0,0)=\sO(1)$, $E(0,1,0,0,0)=\Omega(1)$ and 
$E(1,1,0,0,0)^* = \sT_{\PZ^4}(-2)$.
\end{proof}

\subsection{General Canonical Curves of Genus $6$}

Let $C$ be a general canonical curve of genus $6$. From Mukai's theorem
we obtain a $\PZ^5 \cong \PZ(W) \subset \PZ^9 \cong \PZ(\Lambda_2 V)$ 
and a quadric in $\PZ^5$ such that
\[
    S = \Gr(V,2) \cap \PZ^5 
\]
is a Del Pezzo surface and
\[
    C = S \cap Q
\]

\begin{prop}
On $\PZ^4 \cong \ymin(\Gr(5,2))$ there exists a map of 
vector bundles
\[
       \alpha \colon \sT_{\PZ^4}(-2) \to 6\sO_{\PZ^4}
\]
such that its rank $3$ locus $Z_3(\alpha)$ is the scheme $Z$ of
last scrollar syzygies of $C$. $Z$ is a configuration of
$5$ skew lines in $\PZ^4$.
\end{prop}

\begin{proof} 
$\laststeptext = 1$ so the first scrollar syzygies of $C$ are also the
last scrollar syzygies. The minimal rank first syzygies of $\Gr(5,2)$ are
of rank $4$ and fill the whole space of first syzygies 
$\PZ(\Lambda_{51}^*) = \PZ^4$ as calculated in proposition \ref{yminM6}.

Since $S$ is a general linear section of $\Gr(5,2)$ we can apply 
corollary \ref{determinantal} to obtain a map
\[
         \alpha \colon L \to W\otimes\sO_\ymin
\]
whose rank calculates the rank of syzygies $s \in \ymin$ considered
as syzygies of $S$. In our case this is equal to 
\[
       \alpha \colon \sT_{\PZ^4}(-2) \to 6\sO_{\PZ^4}
\]
since $\PZ(W)=\PZ^6$ and $L = \sT_{\PZ^4}(-2)$. Now $C$ is $S$ intersected 
with a general quadric and therefore $\alpha$
also calculates the rank of $s$ considered a a syzygies of $C$ by 
proposition \ref{generalquadric}.

The last scrollar syzygies of $C$ are first syzygies of rank $3$. 
The argument above shows that the scheme $Z$ of last scrollar
syzygies contains the rank $3$ locus $Z_3(\alpha)$ of $\alpha$.

Since $\PZ(W) \subset \PZ(\Lambda^2 V)$ is a general subspace,
and $L^*=\Omega_{\PZ^4}(2)$ is globally generated, $Z_3(\alpha)$ is 
reduced and of expected dimension 
\[
    \dim Z_3(\alpha) = \dim \PZ^4 - (4-3)(6 -3) = 1. 
\]
On the other hand we are also in the situation of corollary \ref{dimension}
which gives an isomorphism
\[
    \zeta \colon Z \to C^1_4 = \bigcup_{i=1}^{5} \PZ^1.
\]
This shows that $Z_3(\alpha) \subset Z$
is the union of at most $5$ disjoint lines.

Since $Z_3(\alpha)$ is
of expected dimension and we can calculate its class with Porteous formula
\cite{ACGH}[p.86]:
\[
         z_3(\alpha) = \Delta_{6-3,4-3} 
                        \left( 
                              \frac{c_t(6\sO_{\PZ^4})}
                                   {c_t(\sT_{\PZ^4}(-2))}
                        \right) 
                     = \Delta_{3,1} 
                        \left( 
                              \frac{c_t(6\sO_{\PZ^4}(1))}
                                   {c_t(\sT_{\PZ^4}(-1))}
                        \right) 
\]
The Chern polynomials involved are
\[
         c_t(\sT_{\PZ^4}(-1)) = \frac{1}{1-Ht}
\]
as obtained by the Euler-Sequence and 
\[
      c_t(6\sO_{\PZ^4}(1)) = (1+Ht)^6.
\]
This yields
\[
      a = \frac{c_t(6\sO_{\PZ^4}(1))}{c_t(\sT_{\PZ^4}(-1))}
        = (1+Ht)^6(1-Ht) = 1 +5Ht + 9H^2t^2+5H^3t^3 \pm \dots
\]
and
\[
         z_3(\alpha) = \Delta_{3,1}(a) = \det (a_3) = 5H^3.
\]
Since $Z_3(\alpha)$ is reduced this shows that $Z_3(\alpha)$ contains 
all $5$ lines of $Z$. In particular we have $Z$=$Z_3(\alpha)$. 
\end{proof}

\begin{cor}
The ideal sheaf $I_{Z/\PZ^4}$ is resolved by
\begin{align*}
    I_{Z/\PZ^4} &\from 15 \, E(-4,-1,-1,-1,-1) \\ 
            &\from 6 \, E(-5,-1,-1,-1,-2)\\
            &\from E(-6,-1,-1,-1,-3)\\
            &\from 0
\end{align*}
\end{cor}

\begin{proof}
Since $Z=Z_3(\alpha)$ and $\alpha$ drops rank in the expected dimension, the 
ideal sheaf is resolved by the corresponding Eagon-Northcott
complex
\[
    I_{Z/\PZ^4} \from \Lambda^4 W^* \otimes \Lambda^4 L 
        \from \Lambda^5 W^* \otimes \Lambda^4 L \otimes S_1 L
        \from \Lambda^6 W^* \otimes \Lambda^4 L \otimes S_2 L
        \from 0. 
\]
Since $\dim W^* =6$ we have 
\[
    \Lambda^i W^* \otimes \sO = {6 \choose i} \sO.
\]
This gives the above multiplicities.

Furthermore 
\[
    \Lambda^4 L = \Lambda^4 E(1,1,0,0,0)^* 
                  = E(4,1,1,1,1)^* 
                  = E(-4,-1,-1,-1,-1) 
\]
and
\[
    S_i L = S_i E(1,1,0,0,0)^* = E(i,i,0,0,0)^* = E(-i,0,0,0,-i).
\]
Applying these equations to the complex above yields the desired
resolution.
\end{proof}

\begin{thm} \label{main6}
The scheme $Z$ of last scrollar syzygies of a general canonical curve
 $C \subset \PZ^5$ of genus $6$ is a configuration of 
$5$ skew lines in $\PZ^4$ that spans the whole $\PZ^{4}$ of
first syzygies of $C$.
\end{thm}

\begin{proof}
We have to show, that
\[
   Z \subset \PZ(\Lambda_{51}^*)
\]
is non degenerate. It is enough to check $h^0(I_{Z/\PZ^4}(1))=0$.

Since $\sO_{\PZ^4}(1) = E(1,0,0,0,0)$ we obtain
\begin{align*}
   I_{Z/\PZ^4}(1) 
        &\from 15 \, E(-3,-1,-1,-1,-1) \\ 
        &\from 6 \, E(-4,-1,-1,-1,-2) \\
        &\from E(-5,-1,-1,-1,-3) \\
        &\from 0
\end{align*}
We now calculate the cohomology of the above vector bundles using
the theorem of Bott. The fundamental weights of $GL(5)$ are 
$L_1, L_1+L_2, \dots, L_1+\dots+L_5$. The sum of fundamental weights
is therefore $\delta =5L_1 + 4L_2 + 3L_3 + 2L_4 + 1L_5$. We obtain
\[
\begin{array}{|c|l|}\hline
E(\lambda) & \lambda+\delta \\ \hline
E(-3,-1,-1,-1,-1) & (2,3,2,1,0) \\
E(-4,-1,-1,-1,-2) & (1,3,2,1,-1) \\
E(-5,-1,-1,-1,-3) & (0,3,2,1,-2) \\
\hline
\end{array}
\]
For the first two rows we find positive roots $L_i-L_j$ with
$(L_i-L_j,\lambda+\delta)=0$. Therefore all cohomology of these bundles
vanish. For the last row no such root is found, but there
are $3$ roots with $(L_i-L_j,\lambda+\delta)<0$. Therefore the only nonzero
cohomology of $E(-5,-1,-1,-1,-3)$ is $H^3$.

Chasing the diagram
\[
\begin{matrix}
                              & h^0 & h^1 & h^2 & h^3 & h^4 \\
    0 \\
    \uparrow \\
    I_{Z/\PZ^4}(1)                & 0   & *  & 0    & 0   & 0  \\
    \uparrow \\
    15 \, E(-3,-1,-1,-1,-1)   & 0   & 0  & 0    & 0   & 0  \\  
    \uparrow                  & 0   & 0  & *    & 0   & 0  \\ 
    6 \, E(-4,-1,-1,-1,-2)    & 0   & 0  & 0    & 0   & 0  \\ 
    \uparrow \\
    E(-5,-1,-1,-1,-3)         & 0   & 0  & 0    & *   & 0  \\ 
    \uparrow \\
    0 \\
\end{matrix}
\]
we obtain $h^0(I_{Z/\PZ^4}(1))=0$.
\end{proof}

\begin{rem}
Notice that $h^1(I_{Z/\PZ^4}(1))$ does not vanish. 
Therefore $Z$ is not
linearly normal. 
\end{rem}

\section{Genus 7}
\yessubsections

\newcommand{\spinorplus}{{S_{10}^{+}}}
\newcommand{\spinorminus}{{S_{10}^{-}}}
\newcommand{\spinplus}{\spin_5^+}
\newcommand{\spinminus}{\spin_5^-}

\subsection{Syzygies of $M_7$}

Let $V$ be a $10$-dimensional vector space with 
basis $\{ w_1,\dots,w_5,v_1,\dots,v_5\}$. And
$W$ the subspace spanned by the $w_i$.

We consider the Mukai variety
\[
      M_7 = \spinorplus \cong G/P \subset \PZ(\spinplus)
\]
where $G = \Spin(V) = \Spin(10)$ is the $2:1$ spin covering
of $SO(10)$, $\lie{g} = \Lie G = \lie{so_{10}}$ its Lie
algebra and $\spinplus$ a irreducible $16$-dimensional 
spinor representation of $G$. The diagonal matrices 
$D(a_1,a_2,a_3,a_4,a_5,a_1^{-1},a_2^{-1},a_3^{-1},a_4^{-1},a_5^{-1})$
with $a_i \in \CZ$ form a Cartan subalgebra 
$\lie{h} \subset \lie{p} \subset \lie{g}$. The matrices
$H_i = E_{i,i}-E_{i+5,i+5}$ form a basis of $\lie{h}$. Let $L_i$ be
the dual basis of $\lie{h}^*$. Then the positive roots of $G$ are
$L_i \pm L_j$ with $i<j$. The fundamental weight corresponding to
$\spinplus$ is $\frac{1}{2}(L_1+\dots+L_5)$, while the other
fundamental weights are $L_1, L_1+L_2, L_1+L_2+L_3$ 
and $\frac{1}{2}(L_1+\dots+L_4-L_5)$. The corresponding representations
are called $\Lambda_1, \dots, \Lambda_4$ and $\spinminus$ as in \cite{fundamentalKempf}. 

$\spinorplus$ parametrizes one set of $\PZ^4$'s on the smooth
Quadric $Q = v_1w_1 + \dots + v_5w_5$. If
\[
          \pi \colon \Spin(10) \to \SO(10)
\]
is the $2:1$ covering, $\pi(P)$ is the group of orthogonal matrices that
leave a particular $\PZ^4 = \PZ(W)$ invariant, i.e matrices of the form
\[
   \begin{pmatrix}
          *&*&*&*&*&*&*&*&*&*\\ 
          *&*&*&*&*&*&*&*&*&*\\ 
          *&*&*&*&*&*&*&*&*&*\\ 
          *&*&*&*&*&*&*&*&*&*\\ 
          *&*&*&*&*&*&*&*&*&*\\ 
          0&0&0&0&0&*&*&*&*&*\\ 
          0&0&0&0&0&*&*&*&*&*\\ 
          0&0&0&0&0&*&*&*&*&*\\ 
          0&0&0&0&0&*&*&*&*&*\\ 
          0&0&0&0&0&*&*&*&*&*\\ 
   \end{pmatrix}
\]
The semisimple part of $\pi(P)$ is then
\[
 S_{\pi(P)} = \left\{
          \begin{pmatrix}
              A & 0 \\
              0 & (A^{-1})^t
          \end{pmatrix}
          \quad \text{with} \quad
          A \in \GL(W)
          \right\}
\]
         
\begin{prop} 
The  syzygy-numbers of $M_7 = \spinorplus$ are
\[
\begin{matrix}
 1 & -  & -  & -  & -  & - \\
 - & 10 & 16 & -  & -  & - \\
 - & -  & -  & 16 & 10 & - \\
 - & -  & -  & -  & -  & 1\\
\end{matrix}
\]
\end{prop}

\begin{proof} A general linear section $\spinorplus \cap \PZ^6$ 
has the same syzygy numbers as $\spinorplus$ by proposition
\ref{isomorphisms}. The syzygy-numbers of a general canonical curve
$C \subset \PZ^6$ of genus $7$ are the ones claimed above as shown by 
Schreyer in \cite{Sch86}.
\end{proof} 

\begin{prop}
The linear strand of the resolution of $\spinorplus$ is
\[
      I_\spinorplus \from \Lambda_1 \otimes \sO(-2)
                    \from \spinminus \otimes \sO(-3)
\]
\end{prop}

\begin{proof}
From above we have the linear strand
\[
        I_\spinorplus \from V_0 \otimes \sO(-2) 
                      \from V_1 \otimes \sO(-3)
\]
with $\dim V_0 = 10$ and $\dim V_1 = 16$. Now $V_0$ is a
$\Spin(10)$ invariant subset of quadrics in $\PZ(\spinplus)$:
\[
     V_0 \subset S_2(\spinplus) \subset \spinplus \otimes \spinplus
                                = \Lambda_5^+ \oplus \Lambda_3 \oplus \Lambda_1,
\]
where $\Lambda_5^+$ is the irreducible representation
corresponding to the maximal weight vector $L_1+\dots+L_5$. The
representations have dimension $126$,$120$ and $10$ respectively.
Therefore $V_0 = \Lambda_1$. (For the decomposition of the tensor products
see \cite{fundamentalKempf}).

In the next step we know
\[
           V_1 \subset \Lambda_1 \otimes \spinplus 
                       = \lambda_1 \cdot \spinplus \oplus \spinminus
\]
where $\lambda_1 \cdot \spinplus$ denotes the irreducible
representation obtained by adding the maximal weights of $\Lambda_1$
and $\spinplus$. The irreducible summands have dimensions
$144$ and $16$ so that $V_1$ must be equal to $\spinminus$.
\end{proof}

\begin{prop}\label{yminM7}
The scheme $\ymin \subset \PZ(\spinplus)=\PZ^{15}$ 
of minimal rank first syzygies of the Spinor variety $\spinorplus$ 
contains an isomorphic Spinor variety
\[
   \spinorsyz = G/P  
                \subset \PZ(\spinminus)^* 
                \cong \PZ^{15}.
\]
The bundle of linear forms on $\spinorsyz$ is
\[
        L|_\spinorsyz = E(\tfrac{1}{2},\tfrac{1}{2},\tfrac{1}{2},
              \tfrac{1}{2},-\tfrac{1}{2})^*
          = \sB(-1)
\]
where $\sB(-1)$ is the tautological quotient bundle on $\spinorsyz$. 
$\sB(-1)$ is of rank $5$.
\end{prop}

\begin{proof}
From proposition \ref{minorbit} we know, that 
$\ymin \subset \PZ(\spinplus) \cong \PZ^{15}$ must contain the minimal
orbit $G/P$ of $\spin(10)$ in $\PZ(\spinminus)^*$ under the action
\[
       \rho \colon \Spin(10) \to GL((\spinminus)^*).
\] 
This is the spinor variety $\spinorsyz \subset \PZ^{15}$.

To describe the vector bundle of linear forms $L$ on $\spinorsyz$ we have to
determine the action of $P$ on a fiber of $L$. We start by considering
the dual actions $\rho^*$ of $\Spin(10)$, $P$ and $S_{\pi(P)}$
on $\spinminus$.  The lie group of $S_{\pi(P)}$ is
\[
 \Lie S_{\pi(P)} = \left\{
          \begin{pmatrix}
              A & 0 \\
              0 & -A^t
          \end{pmatrix}
          \quad \text{with} \quad
          A \in \lie{gl}(W)
          \right\}
\]
Let $\PZ(W) \subset Q$ be the $\PZ^4$ left invariant by $P$
With this we have the natural 
representations
\[
    \Lambda^{even} W = \Lambda_0 W \oplus \Lambda_2 W \oplus \Lambda_4 W
                     \cong \spinplus
\]
and
\[
    \Lambda^{odd} W =  \Lambda_1 W \oplus \Lambda_3 W \oplus \Lambda_5 W
                     \cong \spinminus
\]
where $\Lie S_P$ acts in the natural way. The 
$w_I = \wedge_{i_I} w_i$ are weight vectors of $G$ and $P$ and
their weights are $\frac{1}{2}(\sum_{i\in I}L_i - \sum_{i \not\in I} L_i)$.
\cite{Fu}[pp. 305-306].

The maximal weight vector
\[
      s = w_1 \wedge w_2 \wedge w_3 \wedge w_4 \wedge w_5 \in \spinminus
\]
is a minimal rank syzygy in $\spinorsyz \subset \ymin$. To
determine the fiber $L|_s$ we restrict the map
\[
       \psi \colon \Lambda_1^* \otimes \sO_{\PZ(\spinminus)^*}(-1)
            \to    \spinplus \otimes  \sO_{\PZ(\spinminus)^*}
\]
from definition \ref{mapofvectorbundles} to $s$. This yields
\[
       \psi|_s = \tildephi(s) \in \Hom(\Lambda_1^*,\spinplus)
                               \cong \Lambda_1 \otimes \spinorplus
\]
where
\[
       \tildephi \colon \spinorminus 
                 \hookrightarrow \Lambda_1 \otimes \spinorplus
\]
is the map defined in \ref{tildephi}. In particular we get
\begin{align*}
      \tildephi(s) 
      = & \tildephi(w_1 \wedge w_2 \wedge w_3 \wedge w_4 \wedge w_5) \\
      = &\quad w_1 \otimes w_2 \wedge w_3 \wedge w_4 \wedge w_5 \\
      &- w_2 \otimes w_1 \wedge w_3 \wedge w_4 \wedge w_5 \\
      &+ w_3 \otimes w_1 \wedge w_2 \wedge w_4 \wedge w_5 \\
      &- w_4 \otimes w_1 \wedge w_2 \wedge w_3 \wedge w_5 \\
      &+ w_5 \otimes w_1 \wedge w_2 \wedge w_3 \wedge w_4
\end{align*}
Consequently the fiber of the line bundle of linear forms is
\begin{align*}
      L|_s  = \Img(\tildephi(s)) =
      \langle
      &w_2 \wedge w_3 \wedge w_4 \wedge w_5,\\ 
      &w_1 \wedge w_3 \wedge w_4 \wedge w_5,\\ 
      &w_1 \wedge w_2 \wedge w_4 \wedge w_5,\\ 
      &w_1 \wedge w_2 \wedge w_3 \wedge w_5,\\ 
      &w_1 \wedge w_2 \wedge w_3 \wedge w_4
      \rangle
\end{align*}
and $L$ is of rank $5$. $S_P$ acts irreducibly on this fiber,
and $w_1 \wedge w_2 \wedge w_3 \wedge w_4$
is the maximal weight vector with weight $\frac{1}{2}(L_1+L_2+L_3+L_4-L_5)$.
Consequently
\[
           E_{\rho^*} = \tfrac{1}{2}E(1,1,1,1,-1)
\]
and
\[      
           L = E_{\rho} = \tfrac{1}{2}E(1,1,1,1,-1)^*
                        = \tfrac{1}{2}E(1,-1,-1,-1,-1).
\]           
Now the tautological sequence on $\spinorplus$ is
\[
          0 \to E(-1,0,0,0,0) 
            \to W \otimes \sO_{\spinorplus}
            \to E(1,0,0,0,0)
            \to 0
\] 
i.e. $\sB=E(1,0,0,0,0)$. Since $\sO(1) =\tfrac{1}{2} E(1,1,1,1,1)$ this
shows that $L|_\spinorsyz=\sB(-1)$.
\end{proof}

\subsection{General Canonical Curves of Genus $7$}

Consider a general canonical curve $C$ of genus $7$. Mukai's 
Theorem provides us with $\PZ^6=\PZ(W) \subset \PZ(\spinplus) \cong \PZ^{15}$ 
(a different $W$ from the last section) such that
\[
         C = \spinorplus \cap \PZ^6
\]
\begin{prop}
On the spinor variety $\spinorsyz \subset \ymin(\spinorplus)$ there 
exists a map of 
vector bundles
\[
       \alpha \colon \sB(-1) \to 7\sO_{\spinorsyz}
\]
such that its rank $3$ locus $Z_3(\alpha)$ is the scheme $Z$ of
last scrollar syzygies of $C$. $Z$ is a ruled surface of degree $84$.
\end{prop}

\begin{proof} 
$\laststeptext = 1$ so the first scrollar syzygies of $C$ are also the
last scrollar syzygies. The minimal rank first syzygies $s \in \spinorsyz$ 
of $\spinorplus$ are
of rank $5$ as calculated in proposition \ref{yminM7}.

Since $C$ is a general linear section of $\spinorplus$ we can apply 
corollary \ref{determinantal} to obtain a map
\[
         \alpha \colon L \to W\otimes\sO_\spinorsyz
\]
whose rank calculates the rank of syzygies $s \in \spinorsyz$ 
considered
as syzygies of $C$. In our case this is equal to 
\[
       \alpha \colon \sB(-1) \to 7\sO_{\spinorsyz}.
\]
The last scrollar syzygies of $C$ are first syzygies of rank $3$. 
The argument above shows that the scheme $Z$ of last scrollar
syzygies contains the rank $3$ locus $Z_3(\alpha)$ of $\alpha$.

Since $\PZ(W) \subset \PZ(\spinplus)$ is a general subspace,
and $L^*=\Omega_{\sB^*}(1)$ is globally generated, $Z_3(\alpha)$ is 
reduced and of expected dimension 
\[
    \dim Z_3(\alpha) = \dim \spinorsyz - (\rank \sB-3)(7 -3) = 2. 
\]
On the other hand we are also in the situation of corollary \ref{dimension}
which gives an isomorphism
\[
    \zeta \colon Z \to C^1_5
\]
with $C^1_5$ a ruled surface over $W^1_5$. Since $C^1_5 \cong Z$ is
irreducible, this shows that $Z_3(\alpha)\subset Z$ is in fact
an equality. In particular $Z$ is
a ruled surface. 

Since $Z=Z_3(\alpha)$ is
of expected dimension and we can calculate its class with Porteous formula
\cite{ACGH}[p.86]:
\[
         z_3(\alpha) = \Delta_{7-3,5-3} 
                        \left( 
                              \frac{c_t(7\sO_{\spinorsyz})}
                                   {c_t(\sB(-1))}
                        \right) 
                     = \Delta_{4,2} 
                        \left( 
                              \frac{c_t(7\sO_{\spinorsyz}(1))}
                                   {c_t(\sB)}
                        \right) 
\]
The cohomology ring of $\spinorsyz$ has been determined by
Ranestad and Schreyer in \cite{Ranestad}[p. 30] as
\[
   H^*(\spinorsyz,\QZ) = \QZ[h,b]/(b^2+8bh^3+8h^6,6h^5b+7h^8)
\]
where $h$ is the class of a hyperplane section and $b$ is the
third Chern class of $\sB^*$. They also give 
the Chern polynomial of $\sB^*$ as  
\[
         c_t(\sB^*) = 1 - 2ht + 2h^2t^2 + bt^3 + (-2h^4-2hb)t^4. 
\]
The Chern polynomials needed for Porteous formula above are
\[
    c_t(\sB) = \frac{1}{c_t(\sB^*)}
\]
as obtained by the tautological sequence and 
\[
      c_t(7\sO_{\PZ^4}(1)) = (1+Ht)^7.
\]
This yields
\begin{align*}
      a &= \frac{c_t(6\sO_{\PZ^4}(1))}{c_t(\sT_{\PZ^4})} 
         = (1+Ht)^7c_t(\sB^*) =\\
        &= 1 +5ht + 9h^2t^2 +(b+7h^3)t^3 + (5bh+5h^4)t^4 
          + (7bh^2+7h^5)t^5 \pm \dots 
\end{align*}
and
\[
         z_3(\alpha) = \Delta_{4,2}(a) 
         = \det 
           \begin{pmatrix}
               a_4 & a_5\\
               a_3 & a_4
           \end{pmatrix}  
         = 7H^8.
\]
Since the degree of the spinor variety is $2g-2=12$ the degree of
$Z=Z_3(\alpha)$ is $7\cdot12=84$. 
\end{proof}

\begin{rem}
On can also calculate the degree of $Z$ via Brill-Noether-Theory.
See \cite{HC} for a formula.
\end{rem}

\begin{cor}
The ideal sheaf $I_{Z/\spinorsyz}$ is resolved by
\begin{align*}
    I_{Z/\spinorsyz} 
        &\from \Lambda_{4}   \otimes S_{1111}  \\
        &\from \Lambda_{5}   \otimes S_{2111}  
             + \Lambda_{41}  \otimes S_{11111}  \\
        &\from \Lambda_{6}   \otimes S_{3111}  
             + \Lambda_{51}  \otimes S_{21111}  \\
        &\from \Lambda_{7}   \otimes S_{4111}  
             + \Lambda_{61}  \otimes S_{31111}  
             + \Lambda_{55}  \otimes S_{22222}  \\
        &\from \Lambda_{71}  \otimes S_{41111}  
             + \Lambda_{65}  \otimes S_{32222}  \\
        &\from \Lambda_{75}  \otimes S_{42222}  
             + \Lambda_{66}  \otimes S_{33222}  \\
        &\from \Lambda_{76}  \otimes S_{43222}  \\
        &\from \Lambda_{77}  \otimes S_{44222}  \\
        &\from 0 
\end{align*}
where $\Lambda_\lambda = \Lambda_\lambda W^*$
and $S_\mu = S_\mu L$.
\end{cor}

\begin{proof}
Since $\alpha \colon L \to 7 \sO_Y$ drops rank in the expected dimension
\[
        \dim Z_3(\alpha) = 10 - (5-3)(7-3) = 2
\] 
the resolution of $I_{Z/\spinorsyz}$ can be calculated using
the methods of Lascoux. In the notation of \cite[Thm 3.3]{Lascoux} 
the above resolution is $k(\alpha,2,0)$ since $\alpha$ drops
rank by two on $Z_3(\alpha)$. 
\end{proof}

This gives

\begin{thm} \label{main7}
The scheme $Z$ of last scrollar syzygies of a general canonical curve
 $C \subset \PZ^6$ of genus $7$ is a linearly normal 
ruled surface of degree $84$
on a spinor variety $\spinorsyz \subset \PZ^{15}$. This ruled
surface spans the whole $\PZ^{15}$ of first syzygies of $C$. 
\end{thm}

\begin{proof}
We have to show, that
\[
      Z \subset \spinorplus \subset \PZ(\spinminus)^*
\]
is non degenerate. It is enough to check $h^0(I_{Z/\spinorsyz}(1))=0$. We tensor
the resolution above by $\sO_Y(1) = \frac{1}{2}E(1,1,1,1,1)$. We calculate
the cohomology of the resulting complex using the theorem of Bott. Notice
that
\begin{align*}
          S_\mu L 
          &= S_\mu \tfrac{1}{2}E(1,-1,-1,-1,-1) \\
          &= S_\mu \bigl( \tfrac{1}{2}E(2,0,0,0,0) \otimes \sO_Y(-1) \bigr) \\
          &= \bigl( S_\mu E(1,0,0,0,0) \bigr) \otimes \sO_Y(-|\mu|) \\
          &= E(\mu) \otimes \sO_Y(-|\mu|) \\
          &= E(\mu-\tfrac{1}{2}|\mu|(1,1,1,1,1)) \\
          &= \tfrac{1}{2}E(2\mu - |\mu|(1,1,1,1,1))
\end{align*}
       
With the sum of fundamental weights 
$\delta = 4\cdot L_1+3\cdot L_2+2 \cdot L_3+1 \cdot L_4+0 \cdot L_5$
we have

\begin{center}
\begin{tabular}{|c|r@{ = }l|c@{$($}r@{,}r@{,}r@{,}r@{,}r|}
\hline
$|\mu|-1$ & \multicolumn{2}{|c|}{Line bundle} & \multicolumn{6}{|c|}{$\delta+\lambda$}  \\
\hline                           
&$\sO(1)$          &$\frac{1}{2}E( 1, 1, 1, 1, 1)$ &&$ 9 $&$ 7 $&$ 5$&$ 3$&$ 1)$ \\ 
 $3$&$S_{1111}L(1) $  &$\frac{1}{2}E(-1,-1,-1,-1,-3)$ &&$ 7 $&$ 5 $&$ 3$&$ 1$&$-3) $\\
 $4$&$S_{2111}L(1) $  &$\frac{1}{2}E( 0,-2,-2,-2,-4)$ &&$ 8 $&$ 4 $&$ 2$&$ 0$&$-4) $\\
 $4$&$S_{11111}L(1) $ &$\frac{1}{2}E(-2,-2,-2,-2,-2)$ &&$ 6 $&$ 4 $&$ 2$&$ 0$&$-2) $\\
 $5$&$S_{3111}L(1) $  &$\frac{1}{2}E( 1,-3,-3,-3,-5)$ &&$ 9 $&$ 3 $&$ 1$&$-1$&$-5) $\\
 $5$&$S_{21111}L(1) $ &$\frac{1}{2}E(-1,-3,-3,-3,-3)$ &&$ 7 $&$ 3 $&$ 1$&$-1$&$-3) $\\
 $6$&$S_{4111}L(1)  $ &$\frac{1}{2}E( 2,-4,-4,-4,-6)$ &&$ 10$&$ 2 $&$ 0$&$-2$&$-6)$ \\
 $6$&$S_{31111}L(1) $ &$\frac{1}{2}E( 0,-4,-4,-4,-4)$ &&$ 8 $&$ 2 $&$ 0$&$-2$&$-4) $\\
 $9$&$S_{22222}L(1) $ &$\frac{1}{2}E(-5,-5,-5,-5,-5)$ &&$ 3 $&$ 1 $&$-1$&$-3$&$-5) $\\
 $7$&$S_{41111}L(1) $ &$\frac{1}{2}E( 1,-5,-5,-5,-5)$ &&$ 9 $&$ 1 $&$-1$&$-3$&$-5)$ \\
$10$&$S_{32222}L(1) $ &$\frac{1}{2}E(-4,-6,-6,-6,-6)$ &&$ 4 $&$ 0 $&$-2$&$-4$&$-6)$ \\
$11$&$S_{42222}L(1) $ &$\frac{1}{2}E(-3,-7,-7,-7,-7)$ &&$ 5 $&$-1 $&$-3$&$-5$&$-7)$ \\
$11$&$S_{33222}L(1) $ &$\frac{1}{2}E(-5,-5,-7,-7,-7)$ &&$ 3 $&$ 1 $&$-3$&$-5$&$-7)$ \\
$12$&$S_{34222}L(1) $ &$\frac{1}{2}E(-4,-6,-8,-8,-8)$ &&$ 4 $&$ 0 $&$-4$&$-6$&$-8)$ \\
$13$&$S_{44222}L(1) $ &$\frac{1}{2}E(-5,-5,-9,-9,-9)$ &&$ 3 $&$ 1 $&$-5$&$-7$&$-9)$ \\ 
\hline
\end{tabular}
\end{center}

where we use the
shorthand notation $(d_1,d_2,d_3,d_4,d_5):=\frac{1}{2}(d_1L_1+\dots +d_5L_5)$ for
the last column.

For the first bundle we have no vanishing root, and all positive
roots $\alpha$ satisfy $(\alpha,\lambda+\delta)>0$. Therefore $i_0=0$ and
all higher cohomology vanishes for $\sO(1)$. For the last bundle there
are $9$ positive roots $\alpha$ with $(\alpha,\lambda+\delta)<0$. Consequently
the only non vanishing cohomology is $H^9(S_{44222}L(1))$. For all remaining
bundles there is at least one integer that appears with opposite signs.
Therefore we have a root $\alpha=L_i+L_j$ with $(\alpha,\lambda+\delta)=0$,
and all cohomology groups  of theses bundles vanish.

Chasing the diagram
\[
\begin{array}{cccccccccccc}
                                         &h^0&h^1&h^2&h^3&h^4&h^5&h^6&h^7&h^8&h^9&h^{10} \\
    0                                   \\
    \uparrow                            \\ 
    I_{Z/\spinorsyz}                              & 0 & 0 & * & 0 & 0 & 0 & 0 & 0 & 0 & 0 & 0\\
    \uparrow                            \\ 
    \Lambda_{4}   \otimes S_{1111}       & 0 & 0 & 0 & 0 & 0 & 0 & 0 & 0 & 0 & 0 & 0\\
    \uparrow                             & 0 & 0 & 0 & * & 0 & 0 & 0 & 0 & 0 & 0 & 0\\ 
    \Lambda_{5}   \otimes S_{2111}  
    + \Lambda_{41}  \otimes S_{11111}    & 0 & 0 & 0 & 0 & 0 & 0 & 0 & 0 & 0 & 0 & 0\\
    \uparrow                             & 0 & 0 & 0 & 0 & * & 0 & 0 & 0 & 0 & 0 & 0\\ 
    \Lambda_{6}   \otimes S_{3111}  
    + \Lambda_{51}  \otimes S_{21111}    & 0 & 0 & 0 & 0 & 0 & 0 & 0 & 0 & 0 & 0 & 0\\
    \uparrow                             & 0 & 0 & 0 & 0 & 0 & * & 0 & 0 & 0 & 0 & 0\\ 
    \Lambda_{7}   \otimes S_{4111}  
    + \Lambda_{61}  \otimes S_{31111}  
    + \Lambda_{55}  \otimes S_{22222}    & 0 & 0 & 0 & 0 & 0 & 0 & 0 & 0 & 0 & 0 & 0\\
    \uparrow                             & 0 & 0 & 0 & 0 & 0 & 0 & * & 0 & 0 & 0 & 0\\ 
    \Lambda_{71}  \otimes S_{41111}  
    + \Lambda_{65}  \otimes S_{32222}    & 0 & 0 & 0 & 0 & 0 & 0 & 0 & 0 & 0 & 0 & 0\\
    \uparrow                             & 0 & 0 & 0 & 0 & 0 & 0 & 0 & * & 0 & 0 & 0\\ 
    \Lambda_{75}  \otimes S_{42222}  
    + \Lambda_{66}  \otimes S_{33222}    & 0 & 0 & 0 & 0 & 0 & 0 & 0 & 0 & 0 & 0 & 0\\
    \uparrow                             & 0 & 0 & 0 & 0 & 0 & 0 & 0 & 0 & * & 0 & 0\\ 
    \Lambda_{76}  \otimes S_{43222}      & 0 & 0 & 0 & 0 & 0 & 0 & 0 & 0 & 0 & 0 & 0\\
    \uparrow                             \\
    \Lambda_{77}  \otimes S_{44222}      & 0 & 0 & 0 & 0 & 0 & 0 & 0 & 0 & 0 & * & 0\\
    \uparrow                            \\
    0 
\end{array}
\]
yields in particular
\[
   h^0(I_{Z/\spinorsyz}(1))=h^1(I_{Z/\spinorsyz}(1))=0.
\]
This shows that the space of scrollar syzygies of $C$
is non degenerate and linearly normal.
\end{proof}

\section{Genus 8}
\yessubsections

This case is very similar to the genus $6$ case, since the Mukai
variety is also a Grassmannian $\Gr(V,2)$ but with $\dim V=6$. 
Let $\{v_1,\dots,v_6\}$ be a basis of $V$.

\subsection{Syzygies of $M_8$}

The Mukai variety for genus $8$ is
\[
          M_8 = \Gr(V,2) = \Gr(6,2) \cong GL(6)/P \subset \PZ(\Lambda_2) \cong \PZ^{14}
\]
The diagonal matrices $\lie{h}$ form a Cartan subalgebra of $\lie{\gl_6}$
and $\lie{p}$. The matrices $H_i= E_{i,i}$ are a basis of $\lie{h}$.
Let $\{L_i\}$ be the dual basis of $\lie{h}^*$.  The positive roots
of $\lie{gl_6}$ are $L_i-L_j$ with $i>j$ and $\omega_i = \sum_{j=0}^i L_j$
are the fundamental weights.

\begin{prop} The syzygy-numbers of $M_8 = \Gr(6,2)$
are
\[
\begin{matrix}
 1 & -  & -  & -  & -  & -  & - \\
 - & 15 & 35 & 21 & -  & -  & - \\
 - & -  & -  & 21 & 35 & 15 & - \\
 - & -  & -  & -  & -  & -  & 1 \\
\end{matrix}
\]
\end{prop}

\begin{proof}  A general linear section $M_8 \cap \PZ^7$ 
has the same syzygy numbers as $M_8$ by proposition
\ref{isomorphisms}. The syzygy-numbers of a general canonical curve
$C \subset \PZ^7$ of genus $8$ are the ones claimed above as shown by 
Schreyer in \cite{Sch86}. 
\end{proof}

\begin{prop}
The linear strand of the resolution of $\Gr(6,2)$ is
\[
   I_{\Gr(6,2)} \from \Lambda_4  \otimes \sO(-2)
              \from \Lambda_{51}  \otimes \sO(-3)
              \from \Lambda_{611}  \otimes \sO(-4)
\]
\end{prop}

\begin{proof}
From above we have a linear strand
\[
   I_{\Gr(6,2)} \from V_0 \otimes \sO(-2)
              \from V_1 \otimes \sO(-3)
              \from V_2 \otimes \sO(-4)
\]
with $\dim V_0 = 15$, $\dim V_1 = 35$ and $\dim V_2 = 21$. 
$V_0$ is a invariant
(not necessarily irreducible)
subspace of quadrics. This gives
\[
    V_0 \subset S_2(\Lambda_2 ) 
        \subset \Lambda_2  \otimes \Lambda_2 
        =  \Lambda_4  \oplus \Lambda_{31}  \oplus \Lambda_{22}  
\]
where the irreducible components have dimensions $15$, $105$ and $105$
respectively. So $V_0 = \Lambda_4 $. Similarly we have
\[
    V_1 \subset \Lambda_4  \otimes \Lambda_2  
         = \Lambda_6 \oplus \Lambda_{51} 
           \oplus \Lambda_{42} 
\]
where the irreducible components have dimensions $1$, $35$ and $175$
respectively. This implies $V_1 = \Lambda_{51}$.

Finally we observe
\[
    V_2 \subset \Lambda_{51} \otimes \Lambda_2
          = \Lambda_{611} \oplus \Lambda_{62} \oplus \Lambda_{521} \oplus \Lambda_{53}
\]
where the irreducible components have dimensions $21$, $15$, $384$ and $105$
respectively. Consequently we have $V_2 = \Lambda_{611}$.
\end{proof}

Using this we get

\begin{prop} \label{yminM8}
The  scheme of minimal rank second syzygies of $\Gr(6,2)$ contains
the minimal orbit
\[
          GL(6)/P 
           \cong \PZ^5 
           \xrightarrow{2-uple} \PZ(\Lambda_{611}^*) 
           \cong \PZ^{20}
\]
The bundle of linear forms on $\PZ^5$ is
\[
        L = E(1,1,0,0,0,0)^* = \sT_{\PZ^5}(-2)
\]
$L$ has rank $5$.
\end{prop}

\begin{proof}
From proposition \ref{minorbit} we know, that 
$\ymin \subset \PZ(\Lambda_{611}^*) \cong \PZ^{20}$ must contain the minimal
orbit of $GL(6)$ in $\PZ(\Lambda_{611}^*)$ under the action
\[
       \rho \colon GL(6) \to GL(\Lambda_{611}^*).
\] 
Here this orbit $GL(5)/P$ is the $2$-uple embedded
$\PZ^5$.

To describe the vector bundle of linear forms on $\PZ^5=GL(6)/P$ we have to
determine the action of $P$ on a fiber of $L$. We start by considering
the dual actions $\rho^*$ of $GL(6)$ and $P$
on $\Lambda_{611}$. The parabolic subgroup in its standard representation
is then the set of
matrices that fix a given $\PZ^0$, i.e matrices of the
form
\[
   \begin{pmatrix}
         * & * & * & * & * & * \\
         0 & * & * & * & * & *\\
         0 & * & * & * & * & *\\
         0 & * & * & * & * & *\\
         0 & * & * & * & * & *\\
         0 & * & * & * & * & * 
   \end{pmatrix}
\]
The semisimple part of $P$ is $S_P = \GL(1) \times \GL(5)$
where $\GL(1)$ acts on $\langle v_1 \rangle$ and $\GL(5)$ acts on
$\langle v_2,v_3,v_4,v_5,v_6 \rangle$ in the standard way.

The maximal weight vector
\[
   s = \begin{tabular}{|c|c|c|}\hline
   1 & 1 & 1 \\ \hline
   2 \\ \cline{1-1}
   3 \\ \cline{1-1}
   4 \\ \cline{1-1}
   5 \\ \cline{1-1}
   6 \\ \cline{1-1}
  \end{tabular} 
   \in \Lambda_{611}
\]
is a syzygy of minimal rank in $\PZ^5 \hookrightarrow \PZ^{20}$.
To determine the fiber $L|_s$ of the bundle of linear forms 
over $s$ we restrict the map
\[
    \psi \colon \Lambda_{51}^* \otimes \sO_{\PZ(\Lambda_{611}^*)}(-1)
                \to \Lambda_2 \otimes \sO_{\PZ(\Lambda_{611}^*)}
\]
from definition \ref{mapofvectorbundles} to $s$. This gives
\[
             \psi|_s =\tildephi(s) \in \Hom(\Lambda_{51}^*,\Lambda_2) 
                                       \cong \Lambda_{51} \otimes \Lambda_2
\]
where
\[
             \tildephi \colon 
             \Lambda_{611} 
             \hookrightarrow 
             \Lambda_{51} \otimes \Lambda_2.
\]
Using Young diagrams we get

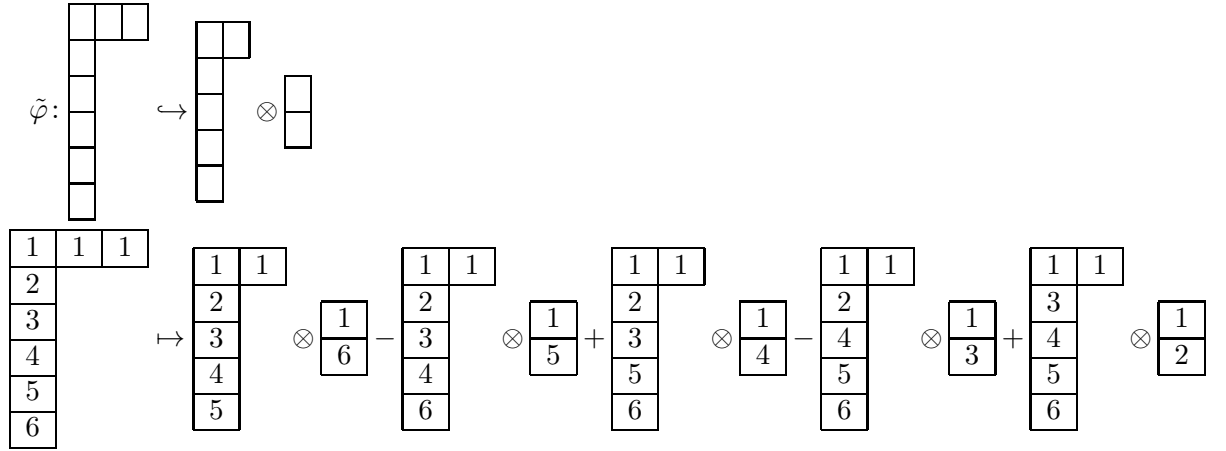
\begin{align*}
   \tildephi \colon \yyyyyy{3}{1}{1}{1}{1}{1}
   &\hookrightarrow  
   \yyyyy{2}{1}{1}{1}{1}\otimes \yy{1}{1}
   \\
   \begin{tabular}{|c|c|c|}\hline
   1 & 1 & 1 \\ \hline
   2 \\ \cline{1-1}
   3 \\ \cline{1-1}
   4 \\ \cline{1-1}
   5 \\ \cline{1-1}
   6 \\ \cline{1-1}
   \end{tabular} 
   &\mapsto
   \begin{tabular}{|c|c|}\hline
   1 & 1\\ \hline
   2 \\ \cline{1-1}
   3 \\ \cline{1-1}
   4 \\ \cline{1-1}
   5 \\ \cline{1-1}
   \end{tabular} 
   \otimes
   \begin{tabular}{|c|}\hline
   1 \\ \hline
   6 \\ \hline
   \end{tabular} 
   -
   \begin{tabular}{|c|c|}\hline
   1 & 1 \\ \hline
   2 \\ \cline{1-1}
   3 \\ \cline{1-1}
   4 \\ \cline{1-1}
   6 \\ \cline{1-1}
   \end{tabular} 
   \otimes
   \begin{tabular}{|c|}\hline
   1 \\ \hline
   5 \\ \hline
   \end{tabular} 
   +
   \begin{tabular}{|c|c|}\hline
   1 & 1 \\ \hline
   2 \\ \cline{1-1}
   3 \\ \cline{1-1}
   5 \\ \cline{1-1}
   6 \\ \cline{1-1}
   \end{tabular} 
   \otimes
   \begin{tabular}{|c|}\hline
   1 \\ \hline
   4 \\ \hline
   \end{tabular} 
   -
   \begin{tabular}{|c|c|}\hline
   1 & 1\\ \hline
   2 \\ \cline{1-1}
   4 \\ \cline{1-1}
   5 \\ \cline{1-1}
   6 \\ \cline{1-1}
   \end{tabular} 
   \otimes
   \begin{tabular}{|c|}\hline
   1 \\ \hline
   3 \\ \hline
   \end{tabular} 
   +
   \begin{tabular}{|c|c|}\hline
   1 & 1 \\ \hline
   3 \\ \cline{1-1}
   4 \\ \cline{1-1}
   5 \\ \cline{1-1}
   6 \\ \cline{1-1}
   \end{tabular} 
   \otimes
   \begin{tabular}{|c|}\hline
   1 \\ \hline
   2 \\ \hline
   \end{tabular} 
\end{align*}

Consequently the fiber of the line bundle of linear forms is
\[
    L|_s = \Img \tildephi(s) = 
   \langle 
   \begin{tabular}{|c|}\hline
   1 \\ \hline
   6 \\ \hline
   \end{tabular}, 
   \begin{tabular}{|c|}\hline
   1 \\ \hline
   5 \\ \hline
   \end{tabular}, 
   \begin{tabular}{|c|}\hline
   1 \\ \hline
   4 \\ \hline
   \end{tabular},
   \begin{tabular}{|c|}\hline
   1 \\ \hline
   3 \\ \hline
   \end{tabular}, 
   \begin{tabular}{|c|}\hline
   1 \\ \hline
   2 \\ \hline
   \end{tabular}
   \rangle
   = \langle
   v_1 \wedge v_6, v_1 \wedge v_5, v_1 \wedge v_4, v_1 \wedge v_3, v_1 \wedge v_2
   \rangle
\]
and $L$ is of rank $5$. $S_P$ acts irreducibly on this fiber, and
$v_1 \wedge v_2$ is the maximal weight vector of weight $L_1 + L_2$.

Therefore
\[
          E_{\rho^*} = E(1,1,0,0,0,0)
\]
and
\[
          L = E_\rho = E(1,1,0,0,0,0)^*
\]
Now the tautological sequence of $\PZ^5$ is 
\[
        0 \to E(0,1,0,0,0,0) 
          \to V \otimes \sO_{\PZ^4} 
          \to E(1,0,0,0,0,0) 
          \to 0
\]
i.e. $E(1,0,0,0,0,0)=\sO(1)$, $E(0,1,0,0,0,0)=\Omega(1)$ and 
$E(1,1,0,0,0,0)^* = \sT_{\PZ^5}(-2)$.
\end{proof}

\subsection{General Canonical Curves of Genus $8$}

Let now $C$ be a general canonical curve of genus $8$. From Mukai's
Theorem we obtain a $\PZ^7 \cong \PZ(W) \subset \PZ(\Lambda_2) \cong \PZ^{14}$ such that
\[
           C = \Gr(6,2) \cap \PZ^7
\]
\begin{prop}
On $\PZ^5 \hookrightarrow \PZ^{20}$ there exists a map of 
vector bundles
\[
       \alpha \colon \sT_{\PZ^5}(-2) \to 8\sO_{\PZ^5}
\]
such that its rank $4$ locus $Z_4(\alpha)$ is the scheme $Z$ of
last scrollar syzygies of $C$. $Z$ is a configuration of
$14$ skew conics on the $2$-uple embedding $\PZ^5 \hookrightarrow \PZ^{20}$.
\end{prop}

\begin{proof} 
$\laststeptext = 2$ so the second scrollar syzygies of $C$ are the 
last scrollar syzygies. The minimal rank second syzygies of $\Gr(6,2)$ are
of rank $5$ and fill at least a $2$-uple embedded 
$\PZ^5 \hookrightarrow \PZ^{20}$ as calculated in proposition \ref{yminM8}.

Since $C$ is a general linear section of $\Gr(6,2)$ we can apply 
corollary \ref{determinantal} to obtain a map
\[
         \alpha \colon L \to W\otimes\sO_\ymin
\]
whose rank calculates the rank of syzygies $s \in \ymin$ considered
as syzygies of $C$. In our case this restricts to 
\[
       \alpha \colon \sT_{\PZ^5}(-2) \to 8\sO_{\PZ^5}
\]
on our $\PZ^5$.

The last scrollar syzygies of $C$ are second syzygies of rank $4$. 
The argument above shows that the scheme $Z$ of last scrollar
syzygies contains the rank $4$ locus $Z_4(\alpha)$ of $\alpha$.
 
Since $\PZ(W) \subset \PZ(\Lambda^2 V)$ is a general subspace,
and $L^*=\Omega_{\PZ^5}(2)$ is globally generated, $Z_4(\alpha)$ is 
reduced and of expected dimension 
\[
    \dim Z_4(\alpha) = \dim \PZ^5 - (5-4)(8 -4) = 1. 
\]
On the other hand we are also in the situation of corollary \ref{dimension}
which gives an isomorphism
\[
    \zeta \colon Z \to C^1_5 = \bigcup_{i=1}^{14} \PZ^1
\]
This shows that $Z_4(\alpha)$
is the union of at most $14$ disjoint $\PZ^1$'s. Each of these $\PZ^1$
is the scheme of second minimal rank syzygies of a scroll. These schemes
are rational normal curves of degree $2$ as calculated in proposition
\ref{syzygyfromscroll}. Since they lie on the $2$-uple embedding of
$\PZ^5$ in $\PZ^{20}$ they are the images of lines in $\PZ^5$.

Since $Z_4(\alpha)$ is
of expected dimension and we can calculate its class with Porteous formula
\cite{ACGH}[p.86]:
\[
         z_4(\alpha) = \Delta_{8-4,5-4} 
                        \left( 
                              \frac{c_t(8\sO_{\PZ^5})}
                                   {c_t(\sT_{\PZ^5}(-2))}
                        \right) 
                     = \Delta_{4,1} 
                        \left( 
                              \frac{c_t(8\sO_{\PZ^5}(1))}
                                   {c_t(\sT_{\PZ^5}(-1))}
                        \right) 
\]
The Chern polynomials involved are
\[
         c_t(\sT_{\PZ^5}(-1)) = \frac{1}{1-Ht}
\]
as obtained by the Euler-Sequence and 
\[
      c_t(8\sO_{\PZ^5}(1)) = (1+Ht)^8.
\]
This yields
\[
      a = \frac{c_t(8\sO_{\PZ^5}(1))}{c_t(\sT_{\PZ^5}(-1))}
        = (1+Ht)^8(1-Ht) = 1 +7Ht + 20H^2t^2+28H^3t^3+14H^4t^4  \pm \dots
\]
and
\[
         z_4(\alpha) = \Delta_{4,1}(a) = \det (a_4) = 14H^3.
\]
Since $Z_4(\alpha)$ is reduced this shows that $Z_4(\alpha)$ contains all 
$14$ conics of $Z \subset \PZ^5 \xrightarrow{2-uple} \PZ^{20}$. In particular 
we have $Z=Z_4(\alpha)$. 
\end{proof}

\begin{cor}
The ideal sheaf $I_{Z/\PZ^5}$ is resolved by
\begin{align*}
    I_{Z/\PZ^5} &\from 56 \, E(-5,-1,-1,-1,-1,-1) \\ 
            &\from 28 \, E(-6,-1,-1,-1,-1,-2) \\
            &\from 8 \,  E(-7,-1,-1,-1,-1,-3) \\
            &\from       E(-8,-1,-1,-1,-1,-4)   \\
            &\from 0
\end{align*}
\end{cor}

\begin{proof}
Since $\alpha$ drops rank in the expected dimension
\[
          1 = 5 - (5-4)(8-4),
\]
the 
ideal sheaf is resolved by the corresponding Eagon-Northcott
complex
\[
    I_{Z/\PZ^5}
        \from \Lambda^5 W^* \otimes \Lambda^5 L 
        \from \Lambda^6 W^* \otimes \Lambda^5 L \otimes S_1 L
        \from \Lambda^7 W^* \otimes \Lambda^5 L \otimes S_2 L
        \from \Lambda^8 W^* \otimes \Lambda^5 L \otimes S_3 L
        \from 0. 
\]
Since $\dim W^* =8$ we have 
\[
    \Lambda^i W^* \otimes \sO = {8 \choose i} \sO.
\]
This gives the above multiplicities.

Furthermore 
\[
    \Lambda^5 L = \Lambda^5 E(1,1,0,0,0,0)^* 
                  = E(5,1,1,1,1,1)^* 
                  = E(-5,-1,-1,-1,-1,-1) 
\]
and
\[
    S_i L = S_i E(1,1,0,0,0,0)^* = E(i,i,0,0,0,0)^* = E(-i,0,0,0,0,-i).
\]
Applying these equations to the complex above yields the desired
resolution.
\end{proof}

\begin{thm} \label{main8}
The scheme $Z$ of last scrollar syzygies of a general canonical curve
 $C \subset \PZ^7$ of genus $8$ is a configuration of 
$14$ skew conics that lie on a $2$-uple embedded 
$\PZ^5 \hookrightarrow \PZ^{20}$. $Z$ spans the whole $\PZ^{20}$ of
second syzygies of $C$.
\end{thm}

\begin{proof}
We have to show, that
\[
    Z \subset \PZ(\Lambda_{611}^*) = \PZ^{20}
\]
is non degenerate. It is enough to check $h^0(I_{Z/\PZ^5}(2))=0$.

Since $\sO_{\PZ^5}(2) = E(2,0,0,0,0,0)$, this yields
\begin{align*}
   I_{Z/\PZ^5}(2) 
            &\from 56 \, E(-3,-1,-1,-1,-1,-1) \\ 
            &\from 28 \, E(-4,-1,-1,-1,-1,-2) \\
            &\from 8 \,  E(-5,-1,-1,-1,-1,-3) \\
            &\from       E(-6,-1,-1,-1,-1,-4)   \\
            &\from 0
\end{align*}
We now calculate the cohomology of the above vector bundles using
the theorem of Bott. The fundamental weights of $GL(6)$ are 
$L_1, L_1+L_2, \dots, L_1+\dots+L_6$. The sum of fundamental weights
is therefore $\delta =6L_1 + 5L_2 + 4L_3 + 3L_4 + 2L_5+L_6$. We obtain
\[
\begin{array}{|c|l|}\hline
E(\lambda) & \lambda+\delta \\ \hline
E(-3,-1,-1,-1,-1,-1) & (3,4,3,2,1,0)\\ 
E(-4,-1,-1,-1,-1,-2) & (2,4,3,2,1,-1)\\
E(-5,-1,-1,-1,-1,-3) & (1,4,3,2,1,-2\\
E(-6,-1,-1,-1,-1,-4) & (0,4,3,2,1,-3)  \\
\hline
\end{array}
\]
For the first three rows we find positive roots $L_i-L_j$ with
$(L_i-L_j,\lambda+\delta)=0$. Therefore all cohomology of these bundles
vanish. For the last row no such root is found, but there
are $4$ roots with $(L_i-L_j,\lambda+\delta)<0$. Therefore the only nonzero
cohomology of $E(-6,-1,-1,-1,-1,-4)$ is $H^4$.

Chasing the diagram
\[
\begin{matrix}
                                 & h^0 & h^1 & h^2 & h^3 & h^4 & h^5 \\
    0 \\
    \uparrow \\
    I_{Z/\PZ^5}(2)         & 0   & *   & 0   & 0   & 0   & 0 \\
    \uparrow \\
    56 \, E(-3,-1,-1,-1,-1,-1)   & 0   & 0   & 0   & 0   & 0   & 0 \\  
    \uparrow                     & 0   & 0   & *   & 0   & 0   & 0 \\ 
    28 \, E(-4,-1,-1,-1,-1,-2)   & 0   & 0   & 0   & 0   & 0   & 0 \\ 
    \uparrow                     & 0   & 0   & 0   & *   & 0   & 0 \\ 
    8 \,  E(-5,-1,-1,-1,-1,-3)   & 0   & 0   & 0   & 0   & 0   & 0 \\
    \uparrow \\    
    E(-5,-1,-1,-1,-3)            & 0   & 0   & 0   & 0   & *   & 0  \\ 
    \uparrow \\
    0 
\end{matrix}
\]
we obtain $h^0(I_{Z/\PZ^5}(2))=0$.
\end{proof}

\begin{rem}
Notice that $h^1(I_{Z/\PZ^5}(2))$ does not vanish. 
Therefore $Z$ is not linearly normal. 
\end{rem}

\end{document}